\documentclass[a4paper, 10pt]{article}
\usepackage[margin=1in]{geometry}
\usepackage{graphicx}		
\usepackage{float}			
\usepackage{color}			
\usepackage{xcolor}			
\usepackage{amsfonts}		
\usepackage{amsmath}		
\usepackage{amssymb}		
\usepackage{bm}				
\usepackage{amsthm}			
\usepackage{esint}			
\usepackage{subcaption}		
\usepackage{listings}		
\usepackage{authblk}		
\usepackage{hyperref}		
\usepackage{calligra}		
\usepackage{enumerate}		
\usepackage{mathtools}		
\usepackage{physics}		
\usepackage{authblk}
\usepackage{comment}

\usepackage{titling}

\usepackage{longtable}

\definecolor{deepred}{RGB}{100, 3, 3}
\definecolor{deepblue}{RGB}{3,3,100}
\definecolor{deepgreen}{RGB}{3,100,3}

\newtheorem{thm}{Theorem}[section]
\newtheorem{prop}[thm]{Proposition}
\newtheorem{lem}[thm]{Lemma}
\newtheorem{ass}[thm]{Assumption}

\theoremstyle{definition}

\newtheorem{remark}[thm]{Remark}

\newenvironment{proving}[1][\proofname]{%
	\begin{proof}[#1]$ $\par\nobreak\ignorespaces
	}{%
	\end{proof}
}

\newcommand{\R}{\mathbb{R}}

\newcommand{\N}{\mathbb{N}}

\newcommand{\Z}{\mathbb{Z}}

\newcommand{\BB}{\mathcal{B}}
\newcommand{\I}{\mathbb{I}}

\newcommand{\inv}{^{-1}}
\newcommand{\eps}{\varepsilon}

\newcommand{\NN}{\mathcal{N}}

\newcommand{\esssup}{{\rm ess sup}}
\newcommand{\essinf}{{\rm ess inf}}

\DeclareMathOperator{\Lip}{Lip}
\DeclareMathOperator{\diam}{diam}
\DeclareMathOperator{\Hell}{Hell.}

\numberwithin{equation}{section}

\allowdisplaybreaks

\newcommand\numberthis{\addtocounter{equation}{1}\tag{\theequation}}

\providecommand{\keywords}[1]
{
  \small	
  \textbf{\textit{Keywords---}} #1
}

\newcommand{\be}{\begin{equation}}
\newcommand{\ee}{\end{equation}}

\newcommand{\beqas}{\begin{eqnarray*}}
\newcommand{\eeqas}{\end{eqnarray*}}

\begin{document}
\title{
Bayesian Inversion of Log-normal Eikonal Equations}
\author{Zhan Fei Yeo and Viet Ha Hoang}
\affil{Division of Mathematical Sciences, \\School of Physical and Mathematical Sciences, \\ Nanyang Technological University\\Singapore 637371}
\date{}
\maketitle
\begin{abstract}
We study the Bayesian inverse problem for inferring the log-normal slowness function of the eikonal equation, given noisy observation data on its solution at a set of spatial points. We contribute rigorous proof on the existence and well-posedness of the problem. We then study approximation of the posterior probability measure by solving the truncated eikonal equation, which contains only a finite number of terms in the Karhunen-Loeve expansion of the slowness function, by the Fast Marching Method. The error of this approximation in the Hellinger metric is deduced in terms of the truncation level of the slowness and the grid size in the Fast Marching Method resolution. It is well known that the plain Markov Chain Monte Carlo procedure for sampling the posterior probability is highly expensive. We develop and justify the convergence of a Multilevel Markov Chain Monte Carlo method. Using the heap sort procedure in solving the forward eikonal equation by the Fast Marching Method, our Multilevel Markov Chain Monte Carlo method achieves a prescribed level of accuracy for approximating the posterior expectation of quantities of interest, requiring only an essentially optimal level of complexity. Numerical examples confirm the theoretical results. 
\end{abstract}
\pagebreak
\section{Introduction}
We consider the inverse problem to infer the slowness function of the eikonal equation, given noisy observations on the solution. We follow the Bayesian approach with the Gaussian prior probability measure, where the slowness function is of the log-normal form. We infer the posterior probability measure, which is the conditional probability of the slowness function given the noisy observations. 

The eikonal equation plays an important role in areas such as seismic tomography and computer vision. From observation information on the solution of the forward eikonal equation, which models the shortest travel time of a wave from a source, the inverse problem infers the slowness function. In Deckelnick et al. \cite{deckelnick2011numerical} and Dunbar and Elliott \cite{dunbar2019binary}, the observation noise is assumed to be deterministic. The best fit approach, which finds the best candidate for the slowness function that minimizes the difference between the forward map and the observation data, is studied. Deckelnick et al. \cite{deckelnick2011numerical} assume that the slowness function is a linear expansion of known functions, and find the coefficients of the expansion; while Dunbar and Elliott \cite{dunbar2019binary} assume a binary slowness. The Bayesian approach to inverse problems (see, e.g., Kaipio and Sommersalo \cite{kaipiosomersalo}, Stuart \cite{stuart2010inverse}) regards the observation noise as a random variable which follows a known distribution. The desired physical property (the slowness in our case) is assumed to belong to a prior probability space. The posterior probability is the solution of the inverse problem. While the deterministic approach may need regularization to be well-posed, Bayesian inverse problems always possess a unique solution as long as the forward observation map is measurable (\cite{stuart2010inverse}). For the eikonal equation, the Bayesian framework is employed in the context of a binary slowness function in Dunbar et al. \cite{dunbar2020reconciling}, where the phase-field penalization and the level set approaches are studied. Chada et al. \cite{chada2020tikhonov} employ the ensemble Kalman inversion approach with Tikhonov regularization. 

In this paper, we consider the case where the slowness function is of the log-normal form. Its natural logarithm follows a Gaussian prior probability distribution, and depends on a countable number of normal random variables, as in the Karhunen-Loeve expansion. Log-normal coefficients are popular in forward and inverse uncertainty quantification as they model the situations where the coefficients (the slowness in our case) are always positive and finite, but  can be arbitrarily close to 0 and arbitrarily large. 
Although Bayesian inverse problems for the eikonal equation are important, it appears to us that a solid theoretical framework has not been properly studied. We contribute in this paper the existence of the posterior probability measure, which constitutes the solution of the Bayesian inverse problem, under the log-normal prior. We rigorously establish a well-posedness result in the form of the local Lipschitzness of the posterior with respect to the observation data in the Hellinger distance.  The slowness is expressed as an expansion of normal random variables as in the Kahunen-Loeve expansion, but to numerically solve the forward eikonal equation to sample the posterior, we need to finitely truncate this expansion. We thus study approximation of the posterior probability measure  by the truncated forward problem, where only a finite number of terms in the expansion of the slowness function is considered. Assuming a decaying rate for the sup norm of the coefficient functions of the expansion, an explicit error estimate in the Hellinger metric for the approximation of the posterior, in terms of the finite number of the chosen expansion terms, is derived. We then approximate the posterior probability by solving the truncated eikonal equation by the Fast Marching Method (FMM)( Sethian \cite{sethian1996fast}). The error estimate for the approximation of the posterior, obtained from the numerical solution of the truncated forward eikonal equation, is deduced. It is the sum of the error of finitely truncating the slowness function and the error of the FMM. The results bear some similarity to those of Bayesian inverse problems for forward elliptic equations under log-Gaussian prior, as considered in, e.g, \cite{stuart2010inverse, HoangSchwabmcmclogn, hoang2020analysis}. However, the approaches for establishing the necessary estimates are different, and may also apply to other Hamilton-Jacobi equations. It uses the minimum principle over Lipschitz paths. Further, the $O(h^{1/2})$ convergence rate of the FMM with the mesh size $h$, as established, e.g., in Deckelnick et al. \cite{deckelnick2011numerical}, only holds when the discretization mesh is smaller than an upper bound, which is realization dependent. As this bound can get arbitrarily small in this case of the Gaussian prior,  for a certain mesh size $h$, the  $O(h^{1/2})$ error estimate may not hold for a set of slowness realizations of positive  prior measure. This is different from the situation of the forward elliptic equations considered in \cite{HoangSchwabmcmclogn, hoang2020analysis}, where the finite element error estimate holds (albeit with a realization dependent multiplying constant) as long as the forward solution is sufficiently regular. The dependence of this upper bound on the realizations needs to be carefully studied.  

As the density of the posterior with respect to the Gaussian prior is known  without the normalizing constant, whose numerical approximations may not be possible, Markov Chain Monte Carlo (MCMC) is usually used to approximate the posterior expectation of quantities of interest. However, in the context of Bayesian inverse problems for forward partial differential equations, this process may be prohibitively expensive. A large number of realizations of the forward equation needs to be solved with equally high levels of accuracy, leading to an enormous level of complexity; see, e.g., Hoang et al. \cite{hoang2013complexity} and Hoang et al. \cite{hoang2020analysis} for some quantitative results. This is also the case for the MCMC sampling procedure of the forward eikonal equation, which is solved by the FMM. Multilevel approaches are well known to reduce substantially the computational complexity of approximating expectations of quantities of interest in both forward and inverse uncertainty quantification, and
have attracted significant interests and contributions (see, e.g., the survey papers \cite{Giles2015} and \cite{Giles2018} and the references therein). For Bayesian inverse problems, we mention exemplarily the references \cite{hoang2013complexity}, \cite{Efendiev2015}, \cite{Dodwell2015}, \cite{Hoel2016}, \cite{Beskos2017}, \cite{hoang2020analysis}. For forward elliptic equations with log-normal coefficients, it is well known that the solution to the forward problem is not uniformly bounded for all the realizations. This leads to the possible non-integrability with respect to the prior of the exponential function of the difference of the mismatch function approximations at two consecutive resolution levels. Without taking this into account, the multilevel approximation may be highly inaccurate, as demonstrated numerically in \cite{hoang2020analysis}. To the best of our knowledge, this issue has only been resolved fully rigorously recently in \cite{hoang2020analysis}. This is exactly the case for the forward eikonal equation with a log-normal slowness, whose solution is not uniformly bounded for all the realizations. The unboundedness of the solution to the forward equation, and of the mismatch function, needs to be carefully considered when constructing the Multilevel Markov Chain Monte Carlo (MLMCMC) algorithm for sampling the posterior measure. 
We develop in this paper the  MLMCMC method for the Bayesian inverse problem for the log-normal eikonal equation. The method achieves an optimal convergence rate, using the FMM with multi resolution levels for the forward eikonal equation. It is based on the MLMCMC method developed for elliptic forward equations with log-normal coefficients in Hoang et al. \cite{hoang2020analysis}. The method is an essential modification  for unbounded log-normal coefficients of the method developed for the uniform prior  in Hoang et al. \cite{hoang2013complexity}. 
For independence sampler, the convergence rate of the method is fully rigorously proved. For pCN sampler, if we assume a spectral gap result similar to that in Hairer et al. \cite{hairer2014spectral}, the method is fully rigorously justified. However, as indicated above, the $O(h^{1/2})$ convergence rate of the FMM only holds when the mesh size $h$ is not more than an upper bound, which is realization dependent, and can be arbitrarily close to $0$. For a fixed mesh size $h$, the set of realizations for which this $O(h^{1/2})$ rate may not hold needs to be carefully considered when proving the convergence  of the MLMCMC method.  We successfully extend the  proof in \cite{hoang2020analysis} to show the optimal MLMCMC convergence rate.  
When the heap sort algorithm is employed in solving the forward eikonal equation by the FMM (see Sethian \cite{sethian1996fast}), we show that the complexity of the MLMCMC method is essentially optimal. To obtain a prescribed accuracy level for approximating the posterior expectation of a quantity of interest, the number of floating point operations required is essentially equal to that (with a possible logarithmic multiplying factor) for solving one realization of the forward eikonal equation by the FMM, for the same level of accuracy. Our numerical examples confirm the theoretically established convergence rate and complexity of the MLMCMC method. 

The paper is organized as follows. 

In the next section, we present the setting up of the Bayesian inverse problem with the log-Gaussian prior, where the slowness function depends on a countable number of mutually independent normal random variables. We prove the existence and the well-posedness of the problem. Section \ref{sec:errorbounding} approximates the posterior probability measure by first truncating the slowness function, taking into account only a finite number of normal random variables, and then by solving the resulting truncated forward eikonal equation by the FMM. We show the error estimate for the approximation in the Hellinger distance, which depends on the finite truncating level of the slowness function and the discretization mesh of the FMM. In Section \ref{sec:mlmcmc}, we develop the MLMCMC method for approximating posterior expectations of quantities of interest.  
As mentioned above, to prove the convergence of the MLMCMC method, we need to modify the proof in \cite{hoang2020analysis}, as the $O(h^{1/2})$ error of the forward solver may not hold for realizations of the slowness in a set of a positive prior measure.  
We present the necessary modifications in Appendix \ref{app:justification}. Numerical examples, which confirm the accurateness, the theoretical error estimates and complexity of the MLMCMC method, are presented in Section \ref{sec:numerical}. We illustrate the theoretical convergence rate of the MLMCMC method in examples where we can use a highly accurate Gauss-Hermite quadrature rule to compute a reference posterior expectation of the quantity of interest. Using the heap sort procedure for the FMM,  the numerical examples show that we achieve essentially optimal computational complexity for the MLMCMC, as theoretically predicted. When the slowness function depends on many random variables, where it is not possible to use a  quadrature rule to compute a highly accurate reference posterior expectation, 
we compare the MLMCMC approximated posterior expectation of the forward solution 
to the true solution corresponding to the reference slowness, from which we obtain the observation data. The numerical results demonstrate the accurateness of our MLMCMC method. The MLMCMC is developed in Section \ref{sec:mlmcmc} when the quantity of interest is the solution to the forward equation. The same procedure applies for computing the posterior expectation of the log-normal slowness. The optimal complexity level is achieved when different levels of approximation of the slowness, corresponding to different levels of truncation of the slowness' expansion, is used.
Although we generate the observation data from a simple reference binary slowness, which is not a priori related to the log-normal form, and we consider a rather arbitrary log-Gaussian prior, the posterior expectation of the slowness obtained from our MLMCMC procedure provides a fairly accurate recovery of the reference binary slowness.

Throughout the paper, by $c$ and $C$, we denote  generic constants whose values can change between different appearances. 

\section{Bayesian inverse problem of the log-normal eikonal equation}
\label{sec:BIP}
\subsection{Setting up of the problem}
\label{subsec:setup}
Let $\Omega \subset \R^d$ be an open and bounded Lipschitz domain. Let $(U,\Theta,\gamma)$ be a probability space. We consider the eikonal equation with a random slowness function. Let the {slowness function}  $s : \overline{\Omega}\times U \rightarrow \R^+$ be  continuous with respect to $x\in \bar\Omega$. Let $x_0 \in \Omega$ be the point {source}. We consider the eikonal equation
under the {Soner boundary condition} 
\begin{align*}
	|\nabla T (x,u)| &= s(x,u), \numberthis \label{eq:eikonal}\\
	T(x_0,u) &= 0, \\
	\nabla T(x,u) \cdot n(x) &\geq 0, \numberthis \label{eq: soner} \\
\end{align*}
where $n$ is the outward unit normal vector. The boundary condition  \eqref{eq: soner} prevents reflection from the boundary back into the domain.  The slowness $s$ and the solution $T$ depend on a parameter $u\in U$.  
Formally, we consider the parametric slowness function of the form
\begin{equation}
\label{eq:lognormalslowness}
s (x, u) := s_* (x) + \exp(\overline{s}(x) + \sum_{i=1}^\infty u_i \psi_i (x)),
\end{equation}
 where $s_*, \bar s$ and $\psi_i$ ($i\in \N$) belong to $L^\infty(\Omega)$, and $u_i\in\R$ for $i\in\N$. The function $s_*(x)$ is assumed to be non-negative, i.e. the case where $s_*(x)$ is indentically zero is possible, so the slowness function $s(x,u)$ can be arbitrarily close to 0 and arbitrarily large. We assume that the random variables $u_i$ in \eqref{eq:lognormalslowness} are mutually independent and are distributed according to the normal distribution ${\cal N}(0,1)$ for all $i$. We define $\R^\N\ni u=(u_1,u_2, \ldots)$. We make the following assumption on the functions $\psi_i$.
\begin{ass}\label{ass:1}
The functions $\psi_i$ satisfies $\|\psi_i\|_{L^\infty(\Omega)}<ci^{-p}$ where $p>1$, $c>0$ are constants. In particular, $\sum_{i=1}^\infty \|\psi_i\|_{L^\infty(D)}$ is finite.
\end{ass}
For conciseness, we denote by $b_i=\|\psi_i\|_{L^\infty(\Omega)}$. The space $U\subset\R^\N$ is defined as
\be
U := \left\{ u \in \R^\N : \sum_{i=1}^\infty \abs{u_i}b_i < \infty\right\}.
\label{eq:U}
\ee
Equipping $\R^\N$ with the product $\sigma$-algebra $\bigotimes_{i=1}^\infty{\cal B}(\R)$ where ${\cal B}(\R)$ is the Borel $\sigma$-algebra on $\R$, we define in $\R^\N$ the probability measure (\cite{bogachev1998Gaussian}, \cite{yamasaki1985measures})
\be
\gamma=\bigotimes_{i=1}^\infty {\cal N}(0,1).
\label{eq:rho}
\ee
With Assumption \ref{ass:1}, $\gamma(U)=1$ (see, e.g., \cite{yamasaki1985measures}, \cite{hoang2020analysis}). We then define the $\sigma$-algebra $\Theta$ as the restriction of $\bigotimes_{i=1}^\infty {\cal B}(\R)$ to $U$, and the prior measure on $U$ as the restriction of the measure defined in \eqref{eq:rho} to $U$, still denoted as $\gamma$. 

We consider the Bayesian inverse problem where noisy observations on the solution $T$ of the eikonal equation \eqref{eq:eikonal} at a finite number of spatial points in $\bar\Omega$ are available.
Let $\{x^i \in \overline{\Omega} : i = 1,2, \dots N \}$ be the set of $N$ sampling points. We consider the {forward operator} $G : U \rightarrow \R^N$ 
\begin{equation}
G(u) := (T(x^1,u), \dots T(x^N,u)).
\label{eq:forwardG}
\end{equation}
Let $\vartheta$ be a random noise which follows the Gaussian distribution ${\cal N}(0,\Sigma)$ in $\R^N$, where the $N\times N$ positive definite covariance matrix $\Sigma$ is known. Given the noisy observation 
\be
\delta=G(u)+\vartheta
\label{eq:delta}
\ee
of the forward functional $G(u)$, our purpose is to determine the posterior probability $\gamma^\delta=P(u|\delta)$ on the measurable space $(U, \Theta)$. 
Let the mismatch function be 
\[
\Phi(u, \delta) := \frac{1}{2} |\delta - G(u)|_\Sigma^2,
\]
where $|\cdot|_\Sigma=|\Sigma^{-1/2}\cdot|$, with $|\cdot|$ denoting the Euclidean norm in $\R^N$. We  show in the next subsection that $\gamma^\delta$ is absolutely continuous with respect to the prior $\gamma$; and that the Radon-Nikodym derivative
\begin{equation}
\frac{d\gamma^\delta}{d \gamma} (u) \propto \exp(- \Phi (u, \delta))
\label{eq:RN}
\end{equation}
holds.
\subsection{Existence and well-posedness of the Bayesian inverse problem}
\label{sect: bipwellposed}
We first recall the variational formulation for the solution of the eikonal equation \eqref{eq:eikonal} (see, e.g., \cite{deckelnick2011numerical}). 
We denote by
	\begin{equation}
	\Xi_{x_0} (x) := \{ \xi \in W^{1, \infty}([0, 1], \overline{\Omega}) \; | \; \xi (0) = x_0, \xi (1) = x \}.
	\end{equation}
The solution of the eikonal equation \eqref{eq:eikonal} is represented by
	\begin{equation}
	\label{eq:fermat}
	T(x,u) = \inf_{\xi \in \Xi_{x_0} (x)} \int_0^1 s(\xi(r),u) |{\xi'(r)}|\; dr.
	\end{equation}
We denote by 
\[
s_{\max} (u) = \esssup_{x\in\bar\Omega} s_*+\exp(\esssup_{x\in\bar\Omega}\bar s+\sum_{i=1}^\infty|u_i|b_i),
\]
and
\[
s_{\min}(u)=\essinf_{x\in\bar\Omega} s_*+\exp(\essinf_{x\in\bar\Omega}\bar s -\sum_{i=1}^\infty|u_i|b_i).
\]
For all $x\in \Omega$, we have
\[
s_{\min}(u)\le s(x,u)\le s_{\max}(u).
\]
We recall the following result on the Lipschitzness of the solution $T$ of \eqref{eq:eikonal}     (see \cite{deckelnick2011numerical}).
\begin{prop}
The solution $T$ of the eikonal equation \eqref{eq:eikonal} is Lipschitz. The Lipschitz constant of $T$ satisfies
\be
Lip(T)\le c\,\esssup_{x\in\bar\Omega}s(x),
\label{eq:LipT}
\ee
where $c$ depends only on the domain $\Omega$. 
\end{prop}
We have the following bound for the minimizer of \eqref{eq:fermat}.
\begin{lem}
	\label{lem: lognormalarclengthbound}
Let $\xi \in\Xi_{x_0} (x)$ be a minimizer of \eqref{eq:fermat}. We have
	\begin{equation}
	\int_0^1 |\xi' (r)| \; dr \leq c\exp(2 \sum_{i=1}^\infty \abs{u_i} b_i),
	\end{equation}
where the constant $c$ only depends on the domain $\Omega$. 
\end{lem}
\begin{proving}
From \eqref{eq:LipT}, as $T(x_0,u)=0$ we have $T(x,u)\le cs_{\max}(u)$. On the other hand, from \eqref{eq:fermat}
\[
s_{\min}(u)\int_0^1|\xi'(r)|dr\le T(x,u)\le cs_{\max}(u).
\]
Thus
\[
\int_0^1|\xi'(r)|dr\le cs_{\max}(u)/s_{\min}(u).
\]
The conclusion then follows. 
\end{proving}
We next show the existence of the solution of the Bayesian inverse problem.
\begin{prop}\label{prop:measurability}
The posterior probability $\gamma^\delta$ is absolutely continuous with respect to the prior probability $\gamma$. Further, the Radon-Nikodym derivative \eqref{eq:RN} holds. 
\end{prop}
\begin{proving}
From Theorem 2.1 in \cite{cotter2009bayesian}, it is sufficient to show that the forward map $G:U\to \R^N$ is measurable.  
%
Let $u, v \in U$. 	Let $\xi_{u}$ and $\xi_{v}$ be the minimizers of \eqref{eq:fermat} for $s(\cdot,u)$ and $s(\cdot,v)$ respectively. From \eqref{eq:fermat}, we have
\[
T(x,u)\le \int_0^1s(\xi_v(t),u)|\xi_v'(t)|dt.
\]
 Using the inequality $|\exp(a)-\exp(b)|\le |a-b|(\exp(a)+\exp(b))$, we have
	\begin{align*}
	&T (x, u) - T(x, v)\\
	&\leq \int_0^1 s(\xi_{v} (t),u) |{\xi_{v}'(t)}| - s(\xi_{v} (t),v) |{\xi_{v}'(t)}| \; dt\\
	&\leq \int_0^1 \exp(\esssup\,\overline{s})\abs{\sum_{i =1}^\infty(u_i - v_i) \psi_i (\xi_{v} (t))}\left[ \exp(\sum_{i =1}^\infty u_i \psi_i(\xi_v (t))) + \exp(\sum_{i =1}^\infty v_i \psi_i(\xi_{v} (t))) \right] |{\xi_{v}'(t)}| \; dt\\
		&\leq c\left(\sum_{i=1}^\infty \abs{u_i - v_i}b_i\right)\left(\exp(\sum_{i=1}^\infty|u_i|b_i)+\exp(\sum_{i=1}^\infty|v_i|b_i)\right)\int_0^1|\xi'_{v}(t)|dt\\
		&\le c\left(\sum_{i=1}^\infty|u_i-v_i|b_i\right)\left(\exp(\sum_{i=1}^\infty|u_i|b_i)+\exp(\sum_{i=1}^\infty|v_i|b_i)\right)\exp(2\sum_{i=1}^\infty |v_i|b_i).
	\end{align*}
By the same procedure, we get a similar bound for $T(x,v)-T(x,u)$. Thus
\be
|T(x,u)-T(x,v)|\le c\left(\sum_{i=1}^\infty|u_i-v_i|b_i\right)\left(\exp(\sum_{i=1}^\infty|u_i|b_i)+\exp(\sum_{i=1}^\infty|v_i|b_i)\right)\left(\exp(2\sum_{i=1}^\infty|u_i|b_i)+\exp(2\sum_{i=1}^\infty|v_i|b_i)\right).
\label{eq:TuTv}
\ee
For $J\in\N$, for $u\in U$, we define
$$
u^J =(u_1, u_2, \dots u_J)\in {\mathbb R}^J. 
$$
We	 now show that for $x\in \bar\Omega$,  $T^J(x, u):=T(x,(u_1,\ldots,u_J,0,0,\ldots))$ is  measurable as a map from $U$ to $\mathbb R$. 
From \eqref{eq:TuTv}, $T^J(x,\cdot)$ regarded as a map from $\R^J$ to $\R$ is continuous. For each $X\in{\cal B}(\R)$, there is a set $X^{-1}\in {\cal B}(\R^J)$ such that the preimage  $(T^J)^{-1}(X)$ is the set of $u\in U$ such that $u^J\in X^{-1}$  which is in the sigma algebra $\Theta$. Thus $T^J$ as a map from $U$ to $\R$ is measurable. 	
From \eqref{eq:TuTv}, 
\[
\lim_{J\to\infty}|T^J(x,u)-T(x,u)|=0.
\]
Thus $T(x,\cdot)$ is the pointwise limit of a sequence of measurable functions, and is thus measurable. Hence, the forward map $G$ is measurable. We get the conclusion. 
\end{proving}

We now show that the Bayesian inverse problem is well-posed. 
The well-posedness proof below follows the proofs for similar results in, e.g. \cite{cotter2009bayesian}, \cite{hoang2020analysis}, for other Bayesian inverse problems. We recall the definition of the Hellinger distance
\be
d_{\Hell} (\gamma^\delta, \gamma^{\delta'})^2 := \frac{1}{2} \int_{U} \left( \sqrt{\frac{d\gamma^\delta}{d\gamma}} - \sqrt{\frac{d \gamma^{\delta'}}{d\gamma}} \right)^2 d\gamma.
\label{eq:dHell}
\ee

\begin{prop}
	\label{prop:wellposed}
	For $|\delta|<r$ and $|\delta'|<r$, we have
	\begin{equation}
d_{\Hell}(\gamma^\delta, \gamma^{\delta'}) \le c(r)|\delta - \delta'|_\Sigma.
	\end{equation}
\end{prop}
\begin{proving}
First we show that the normalizing constant 
\[
Z(\delta)=\int_U\exp(-\Phi(u,\delta))d\gamma(u)
\]
in \eqref{eq:RN} is uniformly bounded from zero for all $\delta$ such that $|\delta|<r$. We note that 
\[
\Phi(u,\delta)\le c(|\delta|_\Sigma^2+|G(u)|_\Sigma^2).
\]
As $T(x,u)\le c\,s_{\max}(u)$,
\be
|G(u)|_\Sigma\le c\exp(\sum_{i=1}^\infty|u_i|b_i).
\label{eq:boundGu}
\ee
	Using Lemma \ref{lem:exponentialbounds},  $\int_U|G(u)|_\Sigma^2d\gamma(u)$ is finite. Thus
\[
\int_U\Phi(u,\delta)\le c(r)
\]
when $|\delta|<r$. Let $M$ be a sufficiently large constant. 
We note that
	\be
P(u \; | \Phi(u, \delta) \leq M) \geq 1 - \frac{C}{M},
\label{eq:c/m}
\ee
so 
\be
Z(\delta) = \int_U \exp(-\Phi(u, \delta)) \; d\gamma(u) \geq \left( 1 - \frac{C}{M} \right) \exp(-M).
\label{eq:boundZdelta}
\ee
Thus $Z(\delta)$ is uniformly bounded from 0 for all $\delta$ such that $|\delta|<r$. 
	From \eqref{eq:dHell}
	\begin{align*}
	2d_{\Hell}(\gamma^{\delta}, \gamma^{\delta'})^2 \leq I_1 + I_2
	\end{align*}
	where:
	$$I_1 := \frac{2}{Z(\delta)} \int_U \left[ \exp(- \frac{1}{2}\Phi(u, \delta)) - \exp(-\frac{1}{2}\Phi(u, \delta')) \right]^2 \; d\gamma(u),$$
	$$I_2 := 2 \abs{Z(\delta)^{-\frac{1}{2}} - Z(\delta')^{-\frac{1}{2}}}^2 \int_U \exp(-\Phi(u, \delta')) \; d\gamma(u).$$

	Using the inequality $|\exp(-a)-\exp(-b)|\le |a-b|$ for $a,b>0$ and the Cauchy-Schwartz inequality, we have 
	\begin{align*}
	\abs{\exp(- \frac{1}{2}\Phi(u, \delta)) - \exp(-\frac{1}{2}\Phi(u, \delta'))} &\leq c |\delta - \delta'|_\Sigma |\delta + \delta' - 2G(u)|_{\Sigma}\\
	&\leq c(r + |G(u)|_\Sigma) |\delta - \delta'|_\Sigma.
	\end{align*}
	Thus
	\begin{align*}
	\int_U \abs{\exp(- \frac{1}{2}\Phi(u, \delta)) - \exp(-\frac{1}{2}\Phi(u, \delta'))}^2 \; d\gamma(u) &\leq c |\delta - \delta'|_\Sigma^2 \int_U(r+ |G(u)|_\Sigma)^2 \; d\gamma(u)\\
	\leq c(r) |\delta - \delta'|_\Sigma^2.
	\end{align*}
As $Z(\delta)$ is uniformly bounded away from 0 when $|\delta|<r$,
	$I_1 \leq c(r) |\delta - \delta'|_\Sigma^2$.  
%
	Note that
	\begin{align*}
	I_2 = 2 \abs{Z(\delta)^{-\frac{1}{2}} - Z(\delta')^{-\frac{1}{2}}}^2 Z(\delta') &\leq 2 \max\{ Z(\delta)^{-3}, Z(\delta')^{-3} \} \abs{Z(\delta) - Z(\delta')}^2 Z(\delta')\\
	&\leq c(r) \abs{Z(\delta) - Z(\delta')}^2.
	\end{align*}
With
	\[
	\abs{Z(\delta) - Z(\delta')} \leq \int_U \abs{\exp(-\Phi(u, \delta)) - \exp(-\Phi(u, \delta'))} \; d\gamma(u)\\
	\leq c(r) |\delta - \delta'|_\Sigma,
	\]
we have $I_2\le c|\delta-\delta'|^2$. We then get the conclusion.
\end{proving}
\begin{remark}
We comment on the similarity of the conditions for the existence, uniqueness and wellposedness of the posterior of the Bayesian inverse problem for the forward log-normal eikonal equation in this paper, and the problems with elliptic/parabolic forward equations with a log-normal coefficient. The key inequality we use to show the existence, uniqueness and well-posedness of the Bayesian inverse problem for the forward eikonal equation is \eqref{eq:TuTv}. For  forward log-normal parabolic  equations, we  have a  similar estimate for the observations (see \cite{HoangQuekSchwabmlmcmcparabolic} proof of Proposition 3.1). Here the derivation of \eqref{eq:TuTv} follows a different approach for the eikonal equation, and may also apply for other Hamilton-Jacobi equations. We compare  the functional in \eqref{eq:fermat} at the minimizing paths for two different realizations of the slowness. The second part of the measurability proof in Proposition \ref{prop:measurability} is similar to that in \cite{HoangQuekSchwabmlmcmcparabolic} Proposition 3.1.   
For forward  elliptic equations with a log-normal coefficient, Hoang and Schwab \cite{HoangSchwabmcmclogn} and Hoang et al. \cite{hoang2020analysis} use the measurability of the forward solution, as a map from the prior probability space $U$ in  \eqref{eq:U} to the Sobolev space $H^1$ (with an appropriate boundary condition), which is  established in \cite{Gittelson} (see also \cite{SchwabGittelson}), to show the measurability of the forward map; but it can be shown in a similar fashion as in the proof of Proposition \ref{prop:measurability}, as being shown for the forward log-normal parabolic equations in \cite{HoangQuekSchwabmlmcmcparabolic}, using an estimate similar to \eqref{eq:TuTv}. In the next section, we use \eqref{eq:TuTv} to derive the approximation error of the posterior with respect to the truncation levels of the slowness, in the Hellinger distance. Our proof is similar to that for the corresponding well-posedness result in \cite{HoangSchwabmcmclogn} and \cite{HoangQuekSchwabmlmcmcparabolic}, where a similar estimate to \eqref{eq:TuTv}, for $v=(u_1,\ldots,u_J,0,0,\ldots)$, is used. The key point in this section is establising \eqref{eq:TuTv} for the forward log-normal eikonal equation. 
We note further that the setting in this paper fits into the general assumptions on Bayesian inverse problems in a measurable prior space considered in Hoang \cite{hoang2012bayesian}. We thus also have a similar local Lipschitzness estimate in the Kullback-Leibler distance, which is larger than the Hellinger distance, namely
\[
d_{\rm KL}(\gamma^\delta,\gamma^{\delta'})\le c(r)|\delta-\delta'|
\]
for $|\delta|<r$ and $|\delta'|<r|$, as shown in \cite{hoang2012bayesian}.
\end{remark}
\section{Approximation of the posterior probability measure}
\label{sec:errorbounding}
In this section, we approximate the posterior probability measure by first approximating the forward eikonal equation by the truncated problem, which only takes into account the first $J$ terms in the expansion \eqref{eq:lognormalslowness} of the slowness function. We then consider the numerical approximation of the resulting truncated eikonal equation by the FMM. 
\subsection{The truncated problem}

For $J \in \N$, we consider the {truncated slowness} function 
\begin{equation} \label{eq:truncatedslowness}
s^J (x, u) := s_* (x) + \exp(\overline{s} (x) + \sum_{i = 1}^{J} u_i \psi_i (x)).
\end{equation}
Let $T^J$ be the unique viscosity solution to the truncated eikonal equation
\begin{equation}
|\nabla T^J (x,u)| = s^J (x, u)\ \mbox{for}\ x\in\bar\Omega
\label{eq:truncatedeikonal}
\end{equation}
with $T^J (x_0,u) = 0$, and the Soner boundary condition $\nabla T^J(x,u)\cdot n(x)\ge 0$ for $x\in\partial\Omega$.
We define the {truncated forward operator} $G^J$ as
 \begin{equation}
G^J (u) := (T^{J}(x^1,u), \dots T^{J} (x^N,u)),
\end{equation} 
where $x^i,\  i= 1,2, \dots N$, are the sample points in \eqref{eq:forwardG}. 
We define the {truncated potential}  $\Phi^J$  as
\begin{equation}
\Phi^J(u, \delta) := \frac{1}{2} |\delta - G^J (u)|_{\Sigma}^2.
\label{eq:PhiJ}
\end{equation}
%
%
%
We define the approximated posterior probability corresponding to the truncated eikonal equation as
\begin{equation}
{d\gamma^{J, \delta}\over d\gamma} (u) \propto \exp( - \Phi^J(u, \delta)).
\end{equation}
%
%
From \cite{deckelnick2011numerical} Theorem 2.2, the solution $T^J$ of the truncated eikonal equation \eqref{eq:truncatedeikonal} is Lipschitz. As $T^J(x_0,u)=0$, 
\[
T^J(x,u)\le cs^J_{\max}(u)\le c\exp(\sum_{i=1}^\infty|u_i|b_i).
\]
We thus have the following bound for the truncated forward map $G^J(u)$
\be
|G^J(u)|\le c\exp(\sum_{i=1}^\infty|u_i|b_i).
\label{eq:boundGJu}
\ee
\begin{remark} As $\Phi^J$ only depends on the finite dimensional vector $(u_1,\ldots,u_J)\in {\mathbb R}^J$, the approximated measure $\gamma^{J,\delta}$ is well-defined. A proof identical to that of Proposition \ref{prop:wellposed} shows that
\[
d_{\rm Hell}(\gamma^{J,\delta},\gamma^{J,\delta'})\le c(r)|\delta-\delta'|,
\]
for $\delta$ and $\delta'$ such that $|\delta|<r$ and $|\delta'|<r$.
\end{remark}
We have the following estimate for the approximated posterior measure $\gamma^{J,\delta}$. 
\begin{lem}
	\label{lem: logNeikonalwholetotruncatedbound}
	Under Assumption \ref{ass:1}, 
	\begin{equation}
	d_{\Hell}(\gamma^{J, \delta}, \gamma^{\delta}) \leq C(r) J^{-q}
	\end{equation}
where $q=p-1$. 
\end{lem}
\begin{proving}
	We use a similar procedure  as in the proof of Proposition \ref{prop:wellposed}. We have
	 $$2d_{\Hell}(\gamma^{J, \delta}, \gamma^{\delta}) \leq I_1 + I_2,$$ 
	 where
	$$I_1:= \frac{2}{Z(\delta)} \int_U \abs{\exp(- \frac{1}{2}\Phi(u, \delta)) - \exp(- \frac{1}{2} \Phi^J (u, \delta))}^2 \; d\gamma(u),$$
	$$I_2 := 2 \abs{Z(\delta)^{-\frac{1}{2}} - Z^{J} (\delta)^{-\frac{1}{2}}} Z^J (\delta).$$
	%
	%
	We have
	
	\begin{align*}
	|\exp(- \frac{1}{2} \Phi(u, \delta)) - \exp(- \frac{1}{2} \Phi^{J} (u, \delta))|^2 &\leq c |2\delta - G(u) - G^J (u)|_\Sigma^2 |G^J(u) - G(u)|_\Sigma^2 \\
	&\leq c (|\delta|_\Sigma^2 + |G(u)|_\Sigma^2 + |G^J (u)|_\Sigma^2)|G^J (u) - G(u)|_\Sigma^2.
	\end{align*} 
%
%
From \eqref{eq:TuTv}, we note that for all $x\in\bar\Omega$
\be
|T(x,u)-T^J(x,u)|\le c\left(\sum_{i> J}|u_i|b_i\right)\exp(3\sum_{i=1}^\infty|u_i|b_i).
\label{eq:truncationrate}
\ee
From \eqref{eq:boundGu} and \eqref{eq:boundGJu}, we have
\begin{align*}
\abs{\exp(- \frac{1}{2} \Phi(u, \delta)) - \exp(- \frac{1}{2} \Phi^{J} (u, \delta))}^2\le c(\delta)\left(\sum_{i>J}|u_i|b_i\right)^2\exp(c\sum_{i=1}^\infty|u_i|b_i).
\end{align*}
Thus
\be
I_1\le cJ^{-2q}.
\label{eq:I1}
\ee
The proof for this inequality uses inequalities \eqref{eq:exp1}, \eqref{eq:exp2} and \eqref{eq:exp3}, and is similar to the proof of the similar inequality in Proposition 4.6 of \cite{HoangSchwabmcmclogn}. 
For $I_2$, we observe that
	$$I_2 := 2\abs{Z(\delta)^{-\frac{1}{2}} - Z^J (\delta)^{-\frac{1}{2}}}^2 \int_U \exp(- \Phi^J (u, \delta)) \; d\gamma(u) \le 2\abs{Z(\delta)^{-\frac{1}{2}} - Z^J (\delta)^{-\frac{1}{2}}}^2. $$
%
From \eqref{eq:boundGJu}, $\int_U\Phi^J(u,\delta)d\gamma(u)$ is uniformly bounded for all $J$. The same proof as for bounding $Z(\delta)$ in the proof of Proposition \ref{prop:wellposed} shows that $Z^J(\delta)$ is uniformly bounded below from 0 for all $J$. Thus	
	\begin{align*}
	I_2 &\leq 2\max\{ Z(\delta)^{-\frac{3}{2}}, Z^J (\delta)^{-\frac{3}{2}} \} \abs{Z(\delta) - Z^J (\delta)}^2\\
	&\leq C(\delta) \abs{Z(\delta) - Z^J (\delta)}^2
	\end{align*}
We note that	
\begin{align*}
	\abs{Z(\delta) - Z^J (\delta)}^2 &= \abs{\int_U \exp(-\Phi(u, \delta)) \; d\gamma(u) - \int_U \exp(- \Phi^J (u, \delta)) \; d\gamma(u)}^2\\
	&\le \int_U \abs{\exp(-\Phi(u, \delta)) - \exp(- \Phi^J (u, \delta))}^2  \; d\gamma(u)\le cJ^{-2q}.
	\end{align*}	
	We then have
	\[ 2d_{\Hell}(\gamma^\delta, \gamma^{J, \delta})^2 \leq C(r) J^{-2q}. \]
\end{proving}
\subsection{Numerical approximation of the truncated forward equation by the FMM}
We approximate the posterior probability measure by numerically solving the truncated forward eikonal equation \eqref{eq:truncatedeikonal} by the Fast Marching Method (\cite{deckelnick2011numerical}, \cite{sethian1996fast}). We assume futher in this section that the domain $\Omega$ satisfies the following properties (see \cite{deckelnick2011numerical}): There is a continuous function $\eta \in C(\overline{\Omega}, \R^d)$, and a positive value $\eps > 0$ such that $\forall x \in \overline{\Omega}$, $B_{\eps \alpha}(x + \alpha\eta(x)) \subset \Omega$ for all $0<\alpha<\eps$ where $B_r(y)$ denotes the open ball centred at $y$ with radius $r$. 

We now describe the FMM for the eikonal equation. 
	Let $h > 0$ denote the mesh size. 	
	Let $\Omega_h := \Omega \cap \Z_h^d$ be the set of internal grid points. 
	Let $\Gamma_h$ be the set of points on $\partial\Omega$ of the form  $x_\alpha + s\sigma e_k \in \partial \Omega$ where  $s \in (0,1]$, $e_k$ for  $k \in \{1, 2, \dots d\}$ is a unit vector in the standard basis of $\R^d$, $x_\alpha \in \Omega_h$ and $\sigma\in\{-1,1\}$. We denote by $G_h=\Omega_h\cup\Gamma_h$. 
%
We consider the discrete eikonal equation : Find $T^{J,h}:G_h\to \R$ such that
\be
T^{J,h}({x_0},u)=0,
\label{eq:dinicond}
\ee
	\begin{equation}
	\label{eq:deikonal}
	\sum_{x_\beta \in N_\alpha} \left[ \left( \frac{T^{J,h}(x_\alpha,u) - T^{J,h}(x_\beta,u)}{h_{\alpha\beta}} \right)^+ \right]^2 = s^J(x_\alpha,u)^2,
	\end{equation}
	where 
\[
G_h \supseteq N_\alpha := \begin{cases}
		\{x_\beta \in G_h  : x_\beta \mbox{ is a neighbour of }x_\alpha\} & \text{ if } x_\alpha \in \Omega_h\\
		\{ x_\beta \in \Omega_h : x_\beta \mbox{ is a neighbour of }x_\alpha\}	& \text{ if } x_\alpha \in \Gamma_h;
	\end{cases}
\]
and $h_{\alpha\beta}=|x_\alpha-x_\beta|$. 
We note the following result (see \cite{deckelnick2011numerical} Lemma 2.3).
\begin{lem}
 Problem \eqref{eq:deikonal} with  condition \eqref{eq:dinicond} has a unique solution which satisfies:
	\label{lem: eikonaldiscretelipschitz}
		 $$\displaystyle \forall x_\alpha \in G_h, T^{J,h}(x_\alpha,u) \geq 0$$
		 and
		\be
		\displaystyle \forall x_\alpha, x_\beta \in G_h, \abs{T^{J,h}(x_\alpha,u) - T^{J,h}(x_\beta,u)} \leq c \max_{x \in \overline{\Omega}} \{s^J(x,u)\} |x_\alpha - x_\beta|,
		 \label{eq:LipTJh}
		 \ee
where the constant $c$ only depends on the domain $\Omega$. 
\end{lem}
The FMM determines the solution $T^{J,h}(\cdot,u) : G_h \rightarrow \R$. 
The algorithm terminates in $O(\abs{G_h} \log \abs{G_h})$ operations when the heap sort procedure is employed (see \cite{sethian1996fast}). 
To determine the convergence rate of the FMM, we assume that the functions in \eqref{eq:lognormalslowness} satisfies $s_*\in W^{1,\infty}(\Omega)$, $\bar s\in W^{1,\infty}(\Omega)$ and $\psi_i\in W^{1,\infty}(\Omega)$ for all $i\in\N$. We denote by $\bar b_i=\|\psi_i\|_{W^{1,\infty}(\Omega)}$. We assume further that
\begin{ass}\label{ass:2} The functions $\psi_i$ in \eqref{eq:lognormalslowness} satisfy $\sum_{i=1}^\infty\bar b_i<\infty$.
\end{ass}
The  set of all $u$ such that $\sum_{i=1}^\infty|u_i|\bar b_i$ is finite has $\gamma$ measure 1. To simplify notation, from now on, we identify the prior space $U$ with this set.  
Deckelnick et al. \cite{deckelnick2011numerical} show that for a fixed slowness function, the convergence rate of the FMM is $O(h^{1/2})$. Examining the proof of \cite{deckelnick2011numerical}, we find that 
\begin{prop}
	\label{prop: fmmsqrtconvergencelogN}
There is a positive constant $c$ such that with 
\be
C(u):=
c\exp(c\sum_{i=1}^\infty(b_i+\bar b_i)|u_i|),
\label{eq:Cu}
\ee
and
\be
h_0(u)={1\over c}\exp(-c\sum_{i=1}^\infty(b_i+\bar b_i)|u_i|),
\label{eq:h0u}
\ee
for all $u \in U$, if $h \in (0, h_0(u)],$ then
	\begin{equation}
	\label{eq:fmmrate}
	\max_{x_\alpha \in G_h} \abs{T^J (x_\alpha,u) - T^{J, h}(x_\alpha,u)} \leq C(u) h^{1/2}.
	\end{equation}
\end{prop}
We then define the approximated forward map
\be
G^{J,h} (u) := (T^{J, h}(x^1,u), \dots, T^{J, h} (x^N,u)),
\ee
where $x^i$ for $i\in \{ 1,2, \dots, N \}$ are the sample points in \eqref{eq:forwardG}, which, for simplicity, we assume to belong to $G_h$. 
To use the approximated forward map $G^{J,h}$ for approximating the posterior measure, we need to establish its measurabillity as a map from $U$ to ${\mathbb R}^N$. We have
\begin{prop}
The forward map $G^{J,h}$,  as map from $(U, \Theta)$ to $({\mathbb R}^N, {\mathcal B}({\mathbb R}^N))$, is measurable.
\end{prop}
\begin{proving}
We use the procedure to establish the existence of a solution of the fast marching method \eqref{eq:deikonal} in Deckelnick et al. \cite{deckelnick2011numerical} (see also \cite{deckelnickinterface}). For $u\in U$, we define the map $Z:G_h\to {\mathbb R}$ by $Z_\alpha=M|x_\alpha-x_{\alpha_0}|$ with $Z_{\alpha_0}=0$; $M$ is a constant to be chosen. For $u\in U$, with $u^J=(u_1,\ldots,u_J)\in {\mathbb R}^J$, we denote by $s^J(x,u^J)=s^J(x,u)$, and $T^{J,h}(x_\alpha, u^J)=T^{J,h}(x_\alpha,u)$. 
Let $B$ be an open ball in ${\mathbb R}^J$. We show that for all $x_\alpha\in G_h$, $T^{J,h}(x_\alpha,u^J)$ is a measurable map from $(B,{\mathcal B}(B)) $ to $(\mathbb R, {\mathcal B}(\mathbb{R}))$. As this holds for all open balls $B$, it implies the measurability of $T^{J,h}(x_\alpha,u^J)$ as a map from $({\mathbb R}^J, {\mathcal B}({\mathbb R}^J))$ to $({\mathbb R}, {\mathcal B}({\mathbb R}))$. 
 As $s^J(x,u^J)$ is 
bounded for all 
$x\in\Omega$ and $u^J\in B$, following \cite{deckelnick2011numerical}, we can choose $M$ sufficiently large such that for all $x_\alpha\in G_h$ and all $u^J\in B$,
\[
\sum_{x_\beta\in{\mathcal N}_\alpha}\left[\left({Z_\alpha-Z_\beta\over h_{\alpha\beta}}\right)^+\right]^2\ge s^J(x,u^J)^2.
\]
 We let $U^0(x_\alpha,u^J)=Z_\alpha$. We define recursively $U^k(x_\alpha,u^J): G_h\to {\mathbb R}_{\ge 0}$ for $k\ge 1$ with $U^k(x_{\alpha_0},u^J)=0$ and
\beqas
U^k(x_\alpha,u^J)=\inf\left\{t\ge 0:\ \sum_{x_\beta\in{\mathcal N}_\alpha}\left[\left({t-U^{k-1}(x_\beta,u^J)\over h_{\alpha\beta}}\right)^+\right]^2\ge s^J(x_\alpha,u^J)^2\right\},
\eeqas
(see \cite{deckelnick2011numerical}). Fix $w^J\in B$. Let $\{v^J_n\}_n\subset B$ be such that  $\lim_{n\to\infty}v^J_n=w^J$. We show by induction that for all $k\ge 0$, $U^k(x_\alpha,u^J)$ is uniformly bounded for all $u^J\in B$, and $\lim_{n\to 0}U^k(x_\alpha,v^J_n)=U^k(x_\alpha,w^J)$ for all $x_\alpha\in G_h$. This holds for $k=0$. 
Following \cite{deckelnickinterface}, for $u^J\in B$, let
\[
\eta(t,u^J)=\sum_{x_\beta\in{\mathcal N}_\alpha}\left[\left({t-U^{k-1}(x_\beta,u^J)\over h_{\alpha\beta}}\right)^+\right]^2.
\]
As $U^{k-1}(x_\beta,u^J)\ge 0$, $\eta(0,u^J)=0$ and $\eta(\cdot,u^J)$ is an increasing function. Thus for all $u^J\in B$, as $\lim_{t\to\infty}\eta(t,u^J)=\infty$, there is a unique value $U^k(x_\alpha,u^J)$ such that 
\be
\eta(U^k(x_\alpha,u^J),u^J)=s^J(x_\alpha,u^J)^2.
\label{eq:Uk}
\ee
 As $U^{k-1}(x_\beta,u^J)$ is uniformly bounded for all $u^J\in B$, and $s^J(x_\alpha,u^J)$ is uniformly bounded for all $u^J\in B$, we deduce that $U^k(x_\alpha,u^J)$ is uniformly bounded for all $u^J\in B$. Thus, from the sequence $\{U^k(x_\alpha,v^J_n)\}_n$, we can extract a convergent subsequence, denote by $\{U^k(x_\alpha,v^J_{n_i})\}_{n_i}$. We denote the limit by $V$. From the induction hypothesis, $\lim_{n_i\to\infty}U^ {k-1}(x_\beta,v^J_{n_i})=U^{k-1}(x_\beta,w^J)$. Thus
\[
\lim_{n_i\to\infty} \sum_{x_\beta\in{\mathcal N}_\alpha}\left[\left({U^k(x_\alpha,v^J_{n_i})-U^{k-1}(x_\beta,v^J_{n_i})\over h_{\alpha\beta}}\right)^+\right]^2= \sum_{x_\beta\in{\mathcal N}_\alpha}\left[\left({V-U^{k-1}(x_\beta,w^J)\over h_{\alpha\beta}}\right)^+\right]^2=\eta(V,w^J)
\]
On the other hand, as
\[
\sum_{x_\beta\in{\mathcal N}_\alpha}\left[\left({U^k(x_\alpha,v^J_{n_i})-U^{k-1}(x_\beta,v^J_{n_i})\over h_{\alpha\beta}}\right)^+\right]^2=\eta(U^k(x_\alpha,v^J_{n_i}),v^J_{n_i})=s^J(x_\alpha,v^J_{n_i})^2,
\]
and $\lim_{n_i\to\infty}s^J(x_\alpha,v^J_{n_i})=s^J(x_\alpha,w^J)$, we have 
\be
\eta(V,w^J)=s^J(x_\alpha,w^J)^2.
\label{eq:V}
\ee
 Thus, from \eqref{eq:Uk} and \eqref{eq:V} and the monotonicity of $\eta(\cdot, w^J)$, we have
$
V=U^k(x_\alpha,w^J).
$ Hence, the whole sequence $\{U^k(x_\alpha,v^J_n)\}_n$ converges to $U^k(x_\alpha,w^J)$. The function $U^k(x_\alpha,u^J)$ is, therefore, continuous so is measurable from $B$ to $\mathbb R$. We have that $\{U^k(x_\alpha,u^J)\}_k$ is monotone with respect to $k$ and converges to $T^{J,h}(x_\alpha,u^J)$ (see \cite{deckelnick2011numerical, deckelnickinterface}). As $T^{J,h}(x_\alpha,u^J)$ is the pointwise limit of a sequence of measurable functions, it is measurable.

Let $X\in {\mathcal B}({\mathbb R}^N)$. There is a set $X^{-1}\in {\mathcal B}({\mathbb R}^J)$ such that the preimage of $X$ of the map $G^{J,h}$ from $U$ to ${\mathbb R}^N$ is the set of all $u=(u_1, u_2, \ldots)\in U$ such that $(u_1, \dots,u_J)\in X^{-1}$. This set is $\Theta$ measurable. Thus $G^{J,h}$ as a map from $(U, \Theta)$ to $({\mathbb R}^N, {\mathcal B}({\mathbb R}^N))$ is measurable.
\end{proving}
We define the $J^{\text{th}}$ term {truncated, discrete Bayesian potential} with mesh size $h$ as
\begin{equation}
\label{eq:PhiJh}
\Phi^{J, h}(u, \delta) := \frac{1}{2} |\delta - G^{J, h} (u)|_{\Sigma}^2.
\end{equation}
The approximated posterior probability $\gamma^{J,h,\delta}$ is defined as
\be
\label{eq:rhoJh}
{d\gamma^{J,h,\delta}\over d\gamma}\propto \exp(- \Phi^{J, h}(u, \delta)),
\ee
where the normalizing constant is
\begin{equation}
Z^{J, h} (\delta) := \int_U \exp(- \Phi^{J, h}(u, \delta)) \; d\gamma(u).
\label{eq:ZJh}
\end{equation}
%
From Lemma \ref{lem: eikonaldiscretelipschitz}, we have
\be
T^{J,h}(x_\alpha,u)\le c\max_{x\in\bar\Omega}s^J(x,u)\le c\exp(c\sum_{i=1}^\infty |u_i|b_i),
\label{eq:boundTJh}
\ee
where the constant $c$ only depends on the domain $\Omega$. Thus
	\begin{equation}
	\label{eq:boundGJhu}
	|{G^{J, h}(u)}|_\Sigma \leq c \exp( \sum_{i = 1}^{J} \abs{u_i} b_i) .
	\end{equation}
A proof similar to that  for $Z(\delta)$ in Proposition \ref{prop:wellposed} shows that $Z^{J,h}(\delta)$ is uniformly bounded below from 0.
%
\begin{remark}
As $G^{J,h}(u)$ is measurable, the approximated posterior probability measure $\gamma^{J,h,\delta}$  in \eqref{eq:rhoJh} is well-defined. Using \eqref{eq:boundGJhu}, a proof identical to that of Proposition \ref{prop:wellposed} shows that $\gamma^{J,h,\delta}$ is locally Lipschitz with respect to $\delta$ in the Hellinger distance.
\end{remark}
We now prove the error bound in the Hellinger distance for the approximation of the posterior measure $\gamma^{J,\delta}$ by $\gamma^{J,h,\delta}$. The proof is different from that in \cite{HoangSchwabmcmclogn} for approximating the posterior probability measures of Bayesian inverse problems with forward log-normal  elliptic equations by finite elements. The theoretical error bound \eqref{eq:fmmrate} of the FMM method only holds when the grid size $h$ is not more than the upper bound $h_0(u)$ in \eqref{eq:h0u}, which is not uniform for all $u\in U$. We have:
\begin{lem}
	\label{lem: logNeikonaltruncatedtotruncatedmeshbound}
	%
	The following estimate holds
	\begin{equation}
	d_{\Hell} (\gamma^{J, h, \delta}, \gamma^{J, \delta}) \leq  c(\delta) h^{1/2}.
	\end{equation}
\end{lem}
\begin{proving}
We have
	\begin{align*}
	2d_{\Hell}(\gamma^{J,h, \delta}, \gamma^{J, \delta})^2 \leq I_1 + I_2 
	\end{align*}
	where:
	$$
I_1 := \frac{2}{Z^J (\delta)} \int_{U} \left[ \exp(- \frac{1}{2} \Phi^J (u, \delta)) - \exp(-\frac{1}{2} \Phi^{J, h} (u, \delta)) \right]^2 \; d\gamma(u),
$$
	$$
I_2 := 2\abs{Z^J (\delta)^{-\frac{1}{2}} - Z^{J, h} (\delta)^{-\frac{1}{2}}}^2 \int_U \exp(- \Phi^{J, h} (u, \delta)) \; d\gamma(u).
$$	
Similarly to the proof of Lemma \ref{lem: logNeikonalwholetotruncatedbound}, we have
	\begin{align}
	\Bigg[ \exp(- \frac{1}{2} \Phi^J (u, \delta)) - &\exp(-\frac{1}{2} \Phi^{J, h} (u, \delta)) \Bigg]^2 
	\leq&  c \left( 2|\delta|_\Sigma^2 + |G^J (u)|_\Sigma^2 + |G^{J, h} (u)|_\Sigma^2 \right) |G^{J, h} (u) - G^J (u)|_\Sigma^2.
\label{eq:esta}
	\end{align}	
Fix a mesh size $h>0$. Let
\[
U':= \{ u \in U : h \le h_0 (u) \}.
\]
From Proposition \ref{prop: fmmsqrtconvergencelogN}, if $u\in U'$, then
		\begin{align*}
		|{G^{J, h} (u) - G^J (u)}|_\Sigma \leq c \max_{x_\alpha \in G_h}  |T(x_\alpha,u) - T^{J, h}(x_\alpha,u)|\le cC(u) h^{1/2}.
\end{align*}
From Lemma \ref{lem:exponentialbounds}, this together with estimates \eqref{eq:boundGJu}, \eqref{eq:boundGJhu} and \eqref{eq:esta} give
\[
\int_{U'}\Bigg[ \exp(- \frac{1}{2} \Phi^J (u, \delta)) - \exp(-\frac{1}{2} \Phi^{J, h} (u, \delta)) \Bigg]^2 d\gamma(u)\le c\left(\int_U C(u)^2d\gamma(u)\right)h
	\leq ch.
\]
From Lemma \ref{lem:exponentialbounds}, we have that 
\[
\int_U {1\over h_0(u)}d\gamma(u)<\infty.
\]
For $u\in U\setminus U'$, $1/h_0(u)>1/h$. Thus there is a constant $c$ such that $\gamma(U\setminus U')<ch$. As $\exp(-\frac12\Phi(u,\delta))\le 1$ and $\exp(-\frac12\Phi^J(u,\delta))\le 1$, 
\[
\int_{U\setminus U'}\Bigg[ \exp(- \frac{1}{2} \Phi^J (u, \delta)) - \exp(-\frac{1}{2} \Phi^{J, h} (u, \delta)) \Bigg]^2 
	\leq ch.
\]
Thus $I_1\le c(r)h$. 
Similarly, we have 
$I_2 \leq c(r) h$.
Therefore,
$$2 d_{\Hell}(\gamma^{J, h, \delta}, \gamma^{J, \delta})^2 \leq c(r) h.$$
This implies the final result. 
\end{proving}
%
From Lemmas \ref{lem: logNeikonalwholetotruncatedbound} and \ref{lem: logNeikonaltruncatedtotruncatedmeshbound}, we have the following approximation
\begin{thm}
	If $|\delta|_\Sigma < r$, then
	\begin{equation}
\label{eq:totalerror}
	d_{\Hell} (\gamma^{\delta}, \gamma^{J, h, \delta}) \leq C(r) (J^{-q} + h^{1/2})
	\end{equation}
where $C(r)$ is a constant depending on $r$. 
\end{thm}
%
\section{Multilevel Markov Chain Monte Carlo}
\label{sec:mlmcmc}
We develop the multilevel Markov Chain Monte Carlo (MLMCMC) approach 
in this section. Let $x^*\in\Omega$. We approximate the posterior expectation of $T(x^*,u)$. For simplicity, we assume that $x^*$ belongs to the approximating grid in the FMM. The method follows from that developed for elliptic equations with log-normal coefficients in {Hoang et al.} \cite{hoang2020analysis}. However, as the convergence rate of the FMM is $O(h^{1/2})$, the number of samples chosen for each resolution level to achieve an optimal convergence rate needs to be adjusted correspondingly. We summarize here the MLMCMC approach.
Justification of the convergence rate is presented in Appendix \ref{app:justification}. 

For the mesh size $h=h_l=O(2^{-l})$ for $l\in\N$ in the FMM, to balance the different sources of errors in \eqref{eq:totalerror}, we choose $J=J_l=O( 2^{l/(2q)})$. For conciseness,  we denote the approximated solution $T^{J_l,h_l}(x_\alpha,u)$ of the truncated eikonal equation with $J_l$ terms in the expansion \eqref{eq:lognormalslowness} and mesh size $h_l$ in the FMM as $T^l(x_\alpha,u)$. The mismatch function $\Phi^{J_l,h_l}$ in \eqref{eq:PhiJh} is denoted as $\Phi^l$; and the approximated posterior probability $\gamma^{J_l,h_l,\delta}$ in \eqref{eq:rhoJh} is denoted as $\gamma^l$. 
Now we consider the MLMCMC for the case of the Gaussian prior.
The MLMCMC estimator $E_L^{MLMCMC}[T(x^*,\cdot)]$ of $\mathbb{E}^{\gamma^\delta}[T(x^*,\cdot)]$ is 
\begin{align*}
&E_{L}^{MLMCMC}(T(x^*,\cdot)) \\
=&\sum_{l=l_0+1}^{L} \sum_{l^{\prime}=l_0+1}^{L^{\prime}(l)}\left[E_{M_{l l^{\prime}}}^{\gamma^{l}}\left[A_{1}^{ll^{\prime}}\right]+E_{M_{l l^{\prime}}}^{\gamma^{l-1}}\left[A_{2}^{l l^{\prime}}\right]+E_{M_{l l^{\prime}}}^{\gamma^{l}}\left[A_{3}^{l}\right] \cdot E_{M_{l l^{\prime}}}^{\gamma^{l-1}}\left[A_{4}^{ll^{\prime}}+A_{8}^{l l^{\prime}}\right]\right. \\
&\left.+E_{M_{l l'}}^{\gamma^{l-1}}\left[A_{5}^{l}\right] \cdot E_{M_{l l^{\prime}}}^{\gamma^{l}}\left[A_{6}^{l l^{\prime}}+A_{7}^{ll^{\prime}}\right]\right] \\
&+\sum_{l=l_0+1}^{L}\left[E_{M_{l l_0}}^{\gamma^{l}}\left[A_{1}^{l l_0}\right]+E_{M_{l l_0}}^{\gamma^{l-1}}\left[A_{2}^{l l_0}\right]+E_{M_{l l_0}}^{\gamma^{l}}\left[A_{3}^{l}\right] \cdot E_{M_{l l_0}}^{\gamma^{l-1}}\left[A_{4}^{l l_0}+A_{8}^{l l_0}\right]\right. \\
&\left.+E_{M_{l l_0}}^{\gamma^{l-1}}\left[A_{5}^{l}\right] \cdot E_{M_{l l_0}}^{\gamma^{l}}\left[A_{6}^{l l_0}+A_{7}^{l l_0}\right]\right] \\
&+\sum_{l^{\prime}=l_0+1}^{L^{\prime}(l_0)} E_{M_{l_0l'}}^{\gamma^{l_0}}\left[T^{l^{\prime}}(x^*,u)-T^{l^{\prime}-1}(x^*,u)\right]+E_{M_{l_0l_0}}^{\gamma^{l_0}}\left[T^{l_0}(x^*,u)\right]
\end{align*}
with $l_0\in\N$ being a starting level and $L>l_0$ being the finest resolution level, where 
\begin{align*}
A_{1}^{ll'} &=\left(1-\exp \left(\Phi^{l}(u ; \delta)-\Phi^{l-1}(u ; \delta)\right)\right) Q(u) I^{l}(u),\\
A_{2}^{ll'} &=\left(\exp \left(\Phi^{l-1}(u ; \delta)-\Phi^{l}(u ; \delta)\right)-1\right) Q(u)\left(1-I^{l}(u)\right),\\
A_{3}^{l} &=\left(\exp \left(\Phi^{l}(u, \delta)-\Phi^{l-1}(u, \delta)\right)-1\right) I^{l}(u), \\
A_{4}^{ll'} &=Q(u) I^{l}(u), \\
A_{5}^{l} &=\left(1-\exp \left(\Phi^{l-1}(u, \delta)-\Phi^{l}(u, \delta)\right)\right)\left(1-I^{l}(u)\right), \\
A_{6}^{ll'} &=\exp \left(\Phi^{l}(u, \delta)-\Phi^{l-1}(u, \delta)\right) Q(u) I^{l}(u),\\
A_{7}^{ll'} &=Q(u)\left(1-I^{l}(u)\right), \\
A_{8}^{ll'} &=\exp \left(\Phi^{l-1}(u, \delta)-\Phi^{l}(u, \delta)\right) Q(u)\left(1-I^{l}(u)\right),
\end{align*}
with $Q(u)=T^{l'}(x^*,u)-T^{l'-1}(x^*,u)$ when $l'\ge l_0+1$ and $Q(u)=T^{l_0}(x^*,u)$ when $l'=l_0$.
Here to handle the unboundedness of the solution of the eikonal equation and the mismatch function $\Phi$, we use the truncation function
\begin{equation*}
    I^{l}(u) = \begin{cases}
    1 \quad \text{ if } \quad \Phi^l(u; \delta) - \Phi^{l-1}(u; \delta) \leq 0, \\
    0 \quad \text{ if } \quad \Phi^l(u; \delta) - \Phi^{l-1}(u; \delta) > 0.
    \end{cases}
\end{equation*}
We denote by $E_{M_{ll'}}^{\gamma^l}$  the MCMC sample average of the Markov chain generated by MCMC sampling procedure with the acceptance probability
\be
\alpha^{l}(u, u') = 1 \wedge \exp(\Phi^{l}(u;\delta)-\Phi^{l}(u'; \delta)), \quad u, u'  \in U,
\label{eq:acceptanceprob}
\ee
for the independence sampler and the pCN sampler (see, e.g., Hairer \cite{hairer2014spectral}).
From \eqref{eq:LipTJh}, there are positive constants $c_1$ and $c_2$ such that for all $u\in U$
\[
\Phi^{J,h}(u;\delta)\le c_1+c_2\exp(2\sum_{i=1}^\infty |u_i|b_i).
\]
Following \cite{hoang2020analysis}, we define the probability measure $\bar\gamma$ on $U$ as
\[
{d\bar\gamma\over d\gamma}\propto \exp(-c_1-c_2\exp(2\sum_{i=1}^\infty |u_i|b_i)).
\]
Let ${\bf E}_L$ be the expectation with respect to the probability space generated by the Markov chains in the MLMCMC sampling procedure with the acceptance probability \eqref{eq:acceptanceprob}, and the initial sample $u^{(0)}$ of each Markov chain being distributed accordingly to $\bar\gamma$. 
With the following sampling choices, 
\begin{align*}
&L^{\prime}(l):=L-l, \quad \text { and } \quad M_{l l^{\prime}}:=2^{L-(l+l^{\prime})} \quad \text { for } \quad l \geq 1, l^{\prime} \geq 1 \\
&M_{l l_0}=M_{l_0 l}=2^{L-l} / L^{2} \quad \text { and } \quad M_{l_0l_0}=2^{L} / L^{4}
\end{align*}
the error estimate is
\begin{equation}
    \mathbf{E}_L[|\mathbb{E}^{\gamma^{\delta}}[T(x^*,\cdot)]-{E}_L^{MLMCMC}[T(x^*,\cdot)]|]\leq C(\delta)L^2 2^{-L/2}.
    \label{eq:mlmcmcerror}
\end{equation}
To reduce the effect of the $L^2$ multiplying factor in \eqref{eq:mlmcmcerror}, we can slightly enlarge the sample size $M_{ll'}$ as
\begin{equation}
M_{ll'}=(l+l')^a2^{L-(l+l')}
\label{eq:a}
\end{equation}
for $a>0$.
\begin{table}[htbp!]
    \centering
    \caption{Total MLMCMC error with different sample sizes}
    \begin{tabular}{|c|c|c|c|c|}
    \hline$a$ & $M_{l l^{\prime}}, l, l^{\prime}>1$ & $M_{l l_0}=M_{l_0 l}$ & $M_{l_0l_0}$ & Total error \\
    \hline 0 & $2^{L-\left(l+l'\right)}$ & $2^{L-l} / L^{2}$ & $2^{L} / L^{4}$ & $O\left(L^{2} 2^{-L/2}\right)$ \\
    \hline 2 & $\left(l+l^{\prime}\right)^{2} 2^{L-(l+l^{\prime})}$ & $2^{L-l}$ & $2^{L} / L^{2}$ & $O\left(L \log L 2^{-L/2}\right)$ \\
    \hline 3 & $\left(l+l^{\prime}\right)^{3} 2^{L-(l+l^{\prime})}$ & $l 2^{L-l}$ & $2^{L} / L$ & $O\left(L^{1 / 2} 2^{-L/2}\right)$ \\
    \hline 4 & $\left(l+l^{\prime}\right)^{4} 2^{L-(l+l^{\prime})}$ & $l^{2} 2^{L-l}$ & $2^{L} /\left(\log L^{2}\right)$ & $O\left(\log L 2^{-L/2}\right)$ \\
    \hline
    \end{tabular}
    \label{table:total error summary}
\end{table}
The error of the MLMCMC sampling procedure is recorded in Table \ref{table:total error summary} for some values of $a$. 
The fully rigorous proof of the error estimates of the MLMCMC approximation for Bayesian inverse problems, for forward elliptic equations with log-normal coefficients, is presented in \cite{hoang2020analysis}. For the forward eikonal equation with a log-normal slowness function, it is necessary to modify the proof in \cite{hoang2020analysis} as the theoretical convergence rate of the FMM method in \eqref{eq:fmmrate} only holds when $h\le h_0(u)$, where $h_0(u)$ is not uniform with respect to $u$. We present the necessary modifications of the proof of \cite{hoang2020analysis} in Appendix \ref{app:justification}. 
\begin{remark}\label{rem:complexity}
Using the heap sort algorithm, the number of operations required for solving the eikonal equation on a grid with $m$ total points is $O(m\log m)$ (Sethian \cite{sethian1996fast}). The total number of operations in the MLMCMC algorithm with the number of samples $M_{ll'}=(l+l')^a2^{L-(l+l')}$ in \eqref{eq:a} is
\[
\lesssim\sum_{l=l_0}^L\sum_{l'=l_0}^{L-l}(l2^{dl}+l'2^{dl'})(l+l')^a2^{L-(l+l')}\lesssim L^{a+1}2^{dL}.
\]
\end{remark}
\begin{remark}
An indentical proof as for showing the uniform lower bound of the normalizing constant $Z(\delta)$ in \eqref{eq:boundZdelta} shows that the normalizing constant $Z^{J,h}(\delta)$ in \eqref{eq:ZJh} is uniformly bounded from zero  for all $J$ and $h$; the lower bound is similar to that in \eqref{eq:boundZdelta}.
As shown in \cite{hoang2020analysis} Appendix B, the constant $C(\delta)$ in \eqref{eq:mlmcmcerror} can be bounded by ${\frak a}(2{\frak a}^2+4{\frak a}-4){\mathbb E}^\gamma[({\mathcal V}^{ll'})^2]$, where $1/{\frak a}$ is the uniform lower bound for $Z^{J,h}(\delta)$ and ${\mathcal V}^{ll'}$ is the function on the right hand side of \eqref{eq:Vll} in the present context. In the argument leading to \eqref{eq:boundZdelta}, when the noise covariance $\Sigma$ is small, the upper bound for the mismatch function $\Phi(u,\delta)$ is large. Thus for the right hand side of
\eqref{eq:c/m} to be positive, the constant $M$ needs to be large. This leads to a small lower bound in \eqref{eq:boundZdelta}, i.e. a large value for $\frak a$. As shown in \cite{hoang2020analysis}, ${\mathbb E}^\gamma[({\mathcal V}^{ll'})^2]$ depends on $\sum_{i=1}^\infty (b_i+\bar b_i+b_i^2+\bar b_i^2)$ where $b_i=\|\psi_i\|_{L^\infty(\Omega)}$ and $\bar b_i=\|\nabla\psi_i\|_{L^\infty(\Omega)}$ as defined in Sections \ref{sec:BIP} and \ref{sec:errorbounding}. Thus when the noice covariance $\Sigma$ is small, and when $\sum_{i=1}^\infty(b_i+\bar b_i+b_i^2+\bar b_i^2)$ is large, $C(\delta)$ is large; to achieve a prescribed level of accuracy for the MLMCMC procedure, the finest mesh level $L$ needs to be small. This is also the case for the plain MCMC sampling procedure where the forward equation is solved with equally high levels of accuracy for all the samples. 
\end{remark}

\section{Numerical Experiments}\label{sec:numerical}
In this section, we present numerical experiments to support the theoretical results in the previous sections. 
%
%
%
%
%

First, we demonstrate numerically the convergence rate of the MLMCMC, with respect to the finest level of FMM discretization. For this purpose, we need a highly accurate reference posterior expectation of the quantity of interest. From \eqref{eq:RN}, we use Gauss-Hermite quadrature to compute the normalizing constant, and the integral of the product of the quantity of interest and $\exp(-\Phi(u;\delta))$ with respect to the Gaussian prior. We thus first consider the case where the slowness depends on one random variable.
We consider a slowness function of the form
\begin{equation}
	s({x}, u) := \exp(u \sin(0.5\pi x_1)\sin(0.5\pi x_2)),
\end{equation}
where ${x}=(x_1,x_2)$ belongs to the domain $\Omega :=(-1,1)^2 \subset \R^2$. The slowness depends on the random variable $u \sim {\mathcal N}(0,1)$. We choose 8 data points in \eqref{eq:forwardG}: $(-1/2,-1)$, $(1/2,-1)$, $(-1/2,1)$, $(1/2,1)$, $(-1, -1/2)$, $(-1,1/2)$, $(1,-1/2)$, $(1,1/2)$. 
To generate the observation data $\delta$, we choose a random realization for $u\in\R$ and solve the forward eikonal equation with mesh density $h=2^{-12}$. We then add a randomly generated realisation of the noise $\vartheta$ which follows the Gaussian distribution $\NN(0,0.1^2I)$, where $I$ is the $8\times 8$ identity covariance matrix. The quantity of interest is $T({x}^*)$ where ${x}^*=(1/2,1/2)$. 
To compute the reference posterior expectation in this case, we use a highly accurate Gauss-Hermite quadrature rule, where the solution of the forward eikonal equation at the quadrature nodes are obtained via the FMM with the fine mesh $h = 2^{-12}$, to approximate the integrals with respect to the Gaussian measure on $\R$. As the FMM converges with the reasonably weak rate $O(h^{1/2})$ when the grid size is $h$, we choose larger values of $a$ in \eqref{eq:a} to obtain good convergence. 
We choose the coarsest level $l_0=2$. 

\begin{figure}[htbp!]
    \centering
    \begin{subfigure}{0.49\linewidth}
    \includegraphics[width=\linewidth]{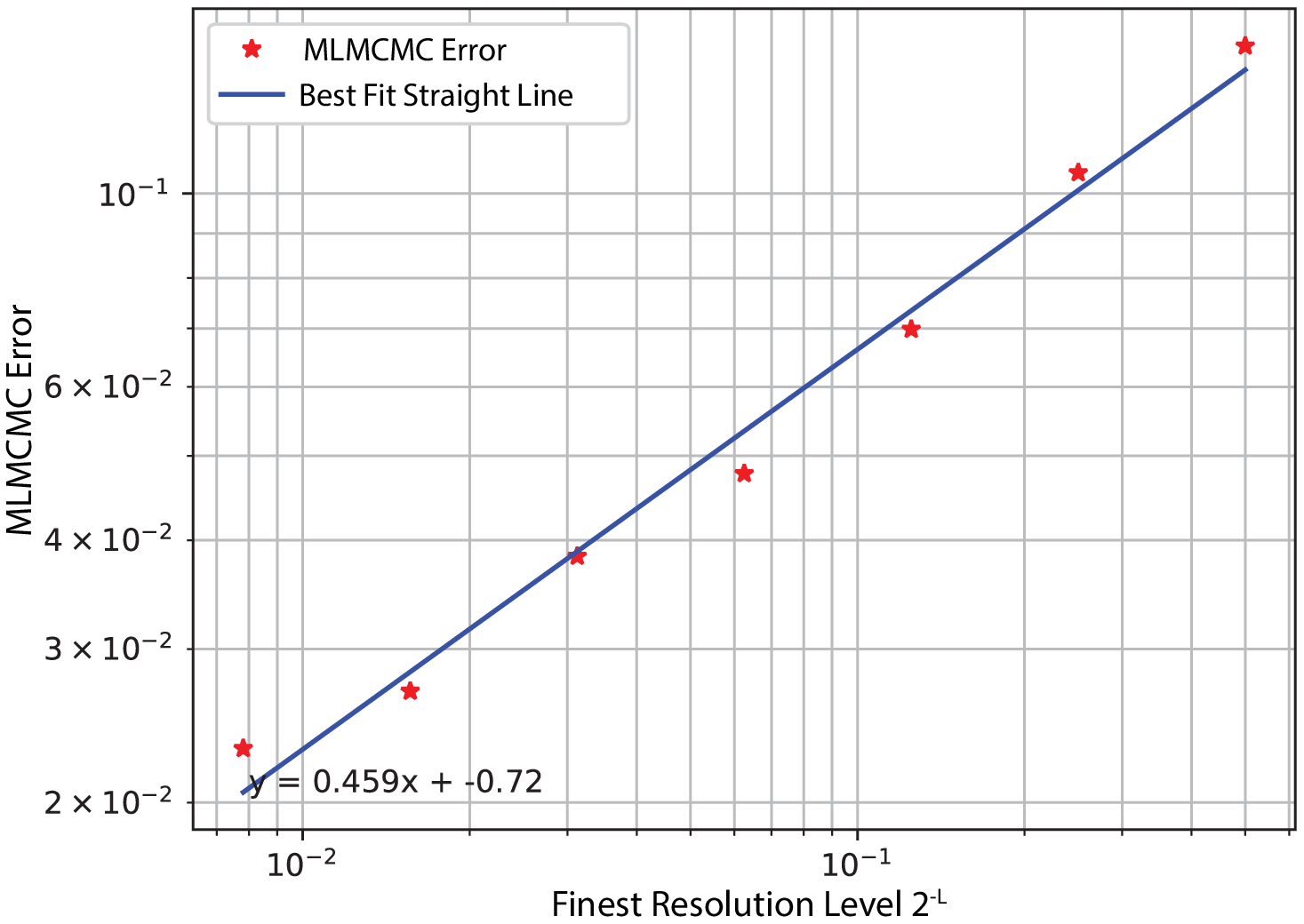}
   	\caption{a = 3}
    \label{fig:alpha=3}
    \end{subfigure}
    \begin{subfigure}{0.49\linewidth}
    \includegraphics[width=0.95\linewidth]{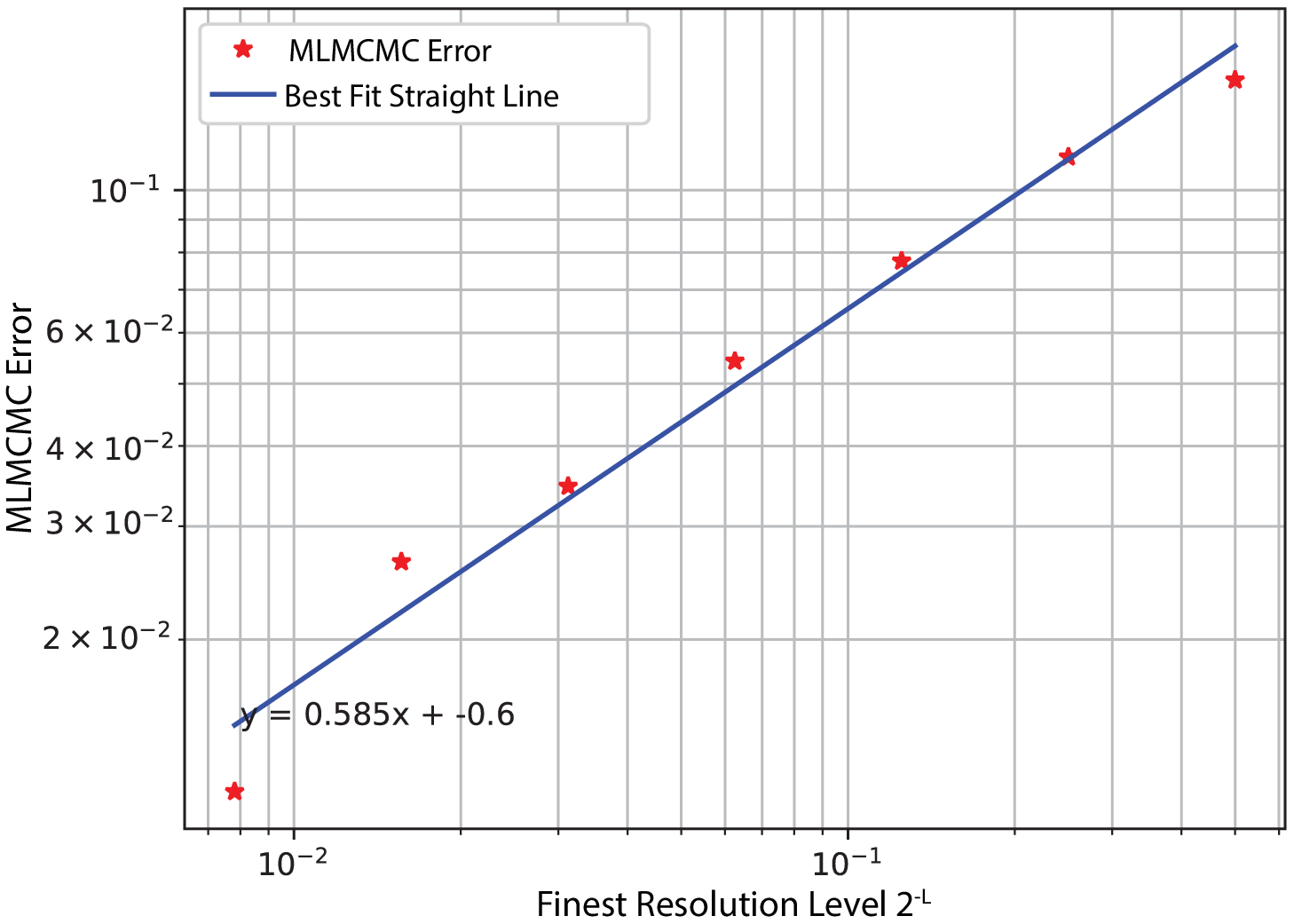}
    \caption{a = 4}
    \label{fig:alpha=4}
    \end{subfigure}
    \caption{Independence Sampler MLMCMC Errors}
    \label{fig:independencemlmcmc}
\end{figure}

Figure \ref{fig:independencemlmcmc} presents the error of the approximated posterior expectation of the quantity of interest ($T((1/2,1/2))$) obtained from the MLMCMC algorithm, where independence sampler is used. The error in the figure is the average of the absolute errors (with respect to the highly accurate reference posterior expectation obtained from the Gauss-Hemite quadrature mentioned above) of 32 indepedent runs of the MLMCMC algorithm. We find that the slope of the best fit straight lines for the value $a=3$ and $a=4$ in \eqref{eq:a} are 0.459 and 0.585 respectively, which are  in reasonable agreement with the theoretically established convergence rate.

Next, we present the error of the MLMCMC where the pCN sampler is used, where the proposal is
\[
v^{(k)}=\sqrt{1-\beta^2}u^{(k)}+\beta \xi,
\]
where $\xi\sim\NN(0,1)$, for generating a proposal from the current sample $u^{(k)}$ of the Markov chain. 

\begin{figure}[htbp!]
    \centering
    \begin{subfigure}{0.49\linewidth}
    \includegraphics[width=\linewidth]{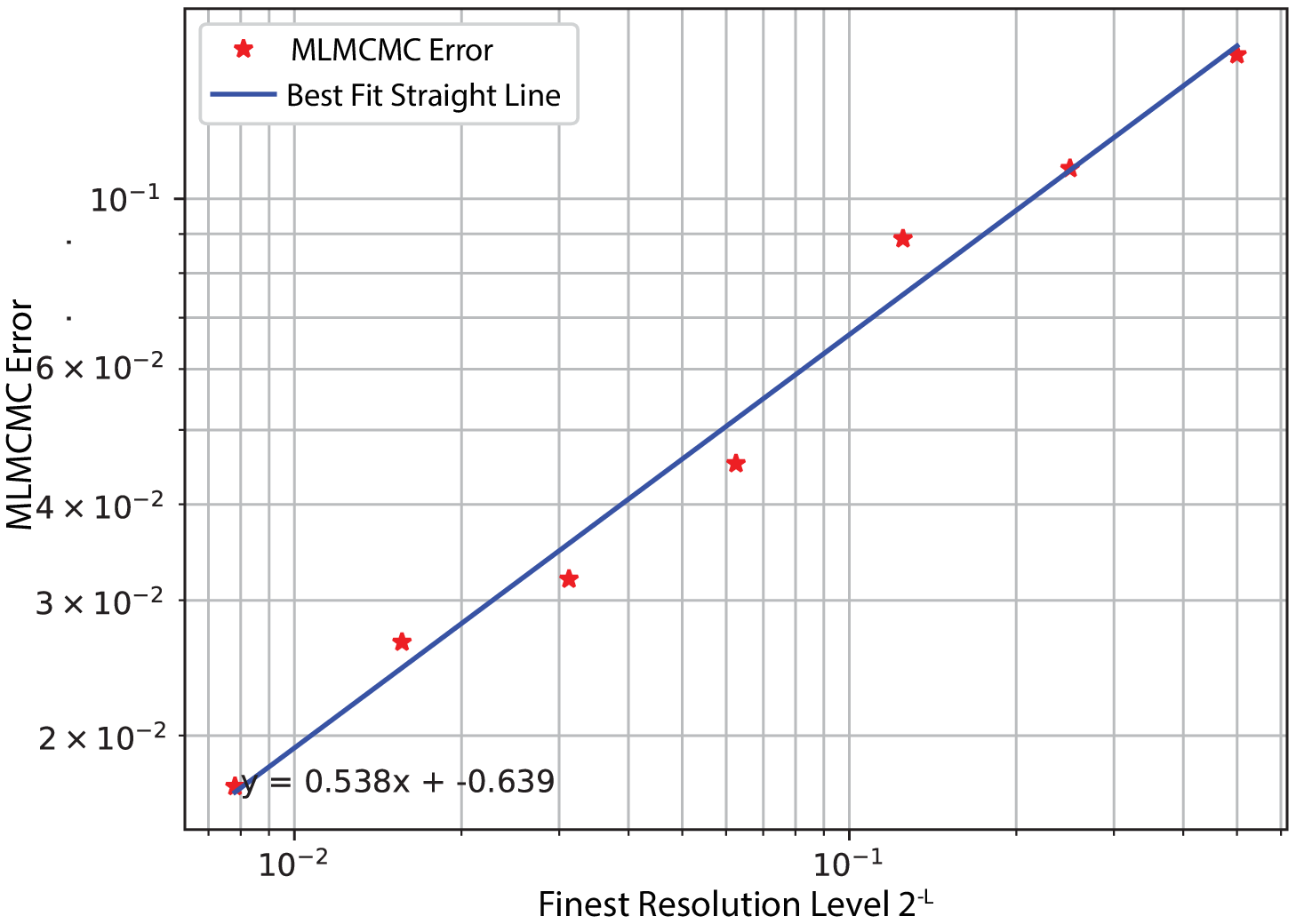}
   	\caption{a = 3}
    \label{fig:pcn08alpha=3}
    \end{subfigure}
    \begin{subfigure}{0.49\linewidth}
    \includegraphics[width=0.95\linewidth]{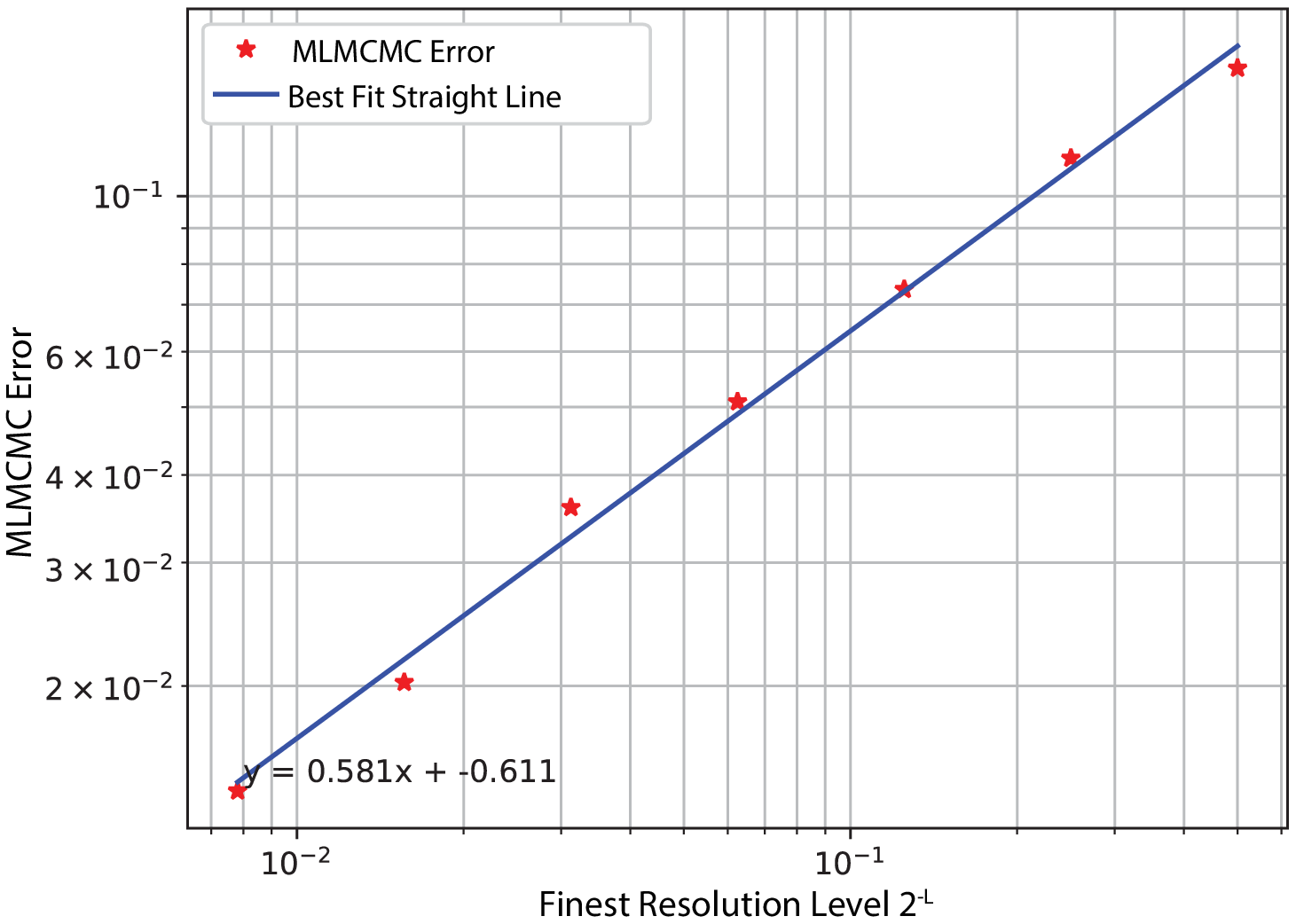}
    \caption{a = 4}
    \label{fig:pcn08alpha=4}
    \end{subfigure}
    \caption{pCN MLMCMC Errors for $\beta = 0.8$}
    \label{fig:pcn08mlmcmc}
\end{figure}

\begin{figure}[htbp!]
    \centering
    \begin{subfigure}{0.465\linewidth}
    \includegraphics[width=\linewidth]{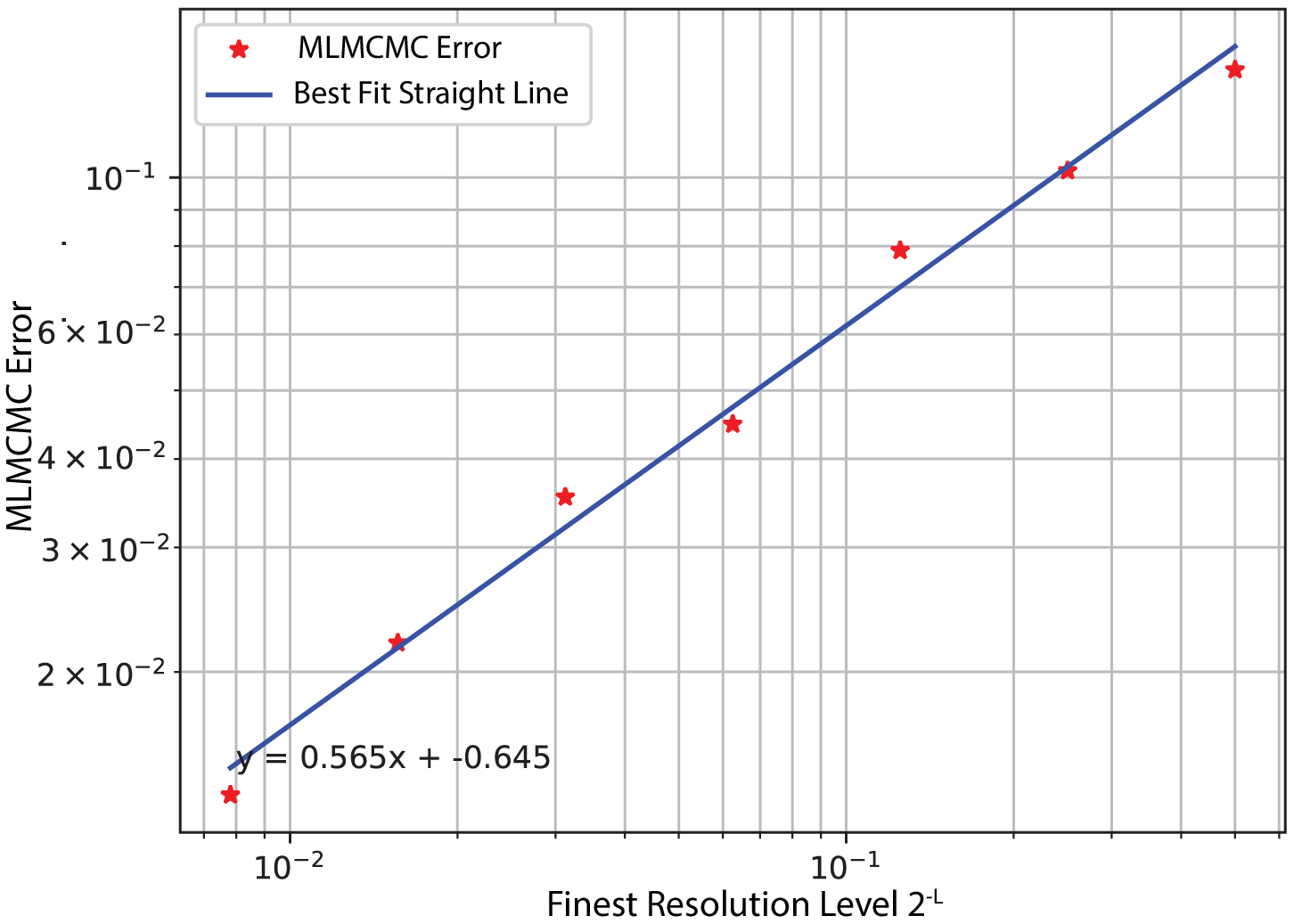}
   	\caption{a = 3}
    \label{fig:pcn05alpha=3}
    \end{subfigure}
    \begin{subfigure}{0.49\linewidth}
    \includegraphics[width=0.95\linewidth]{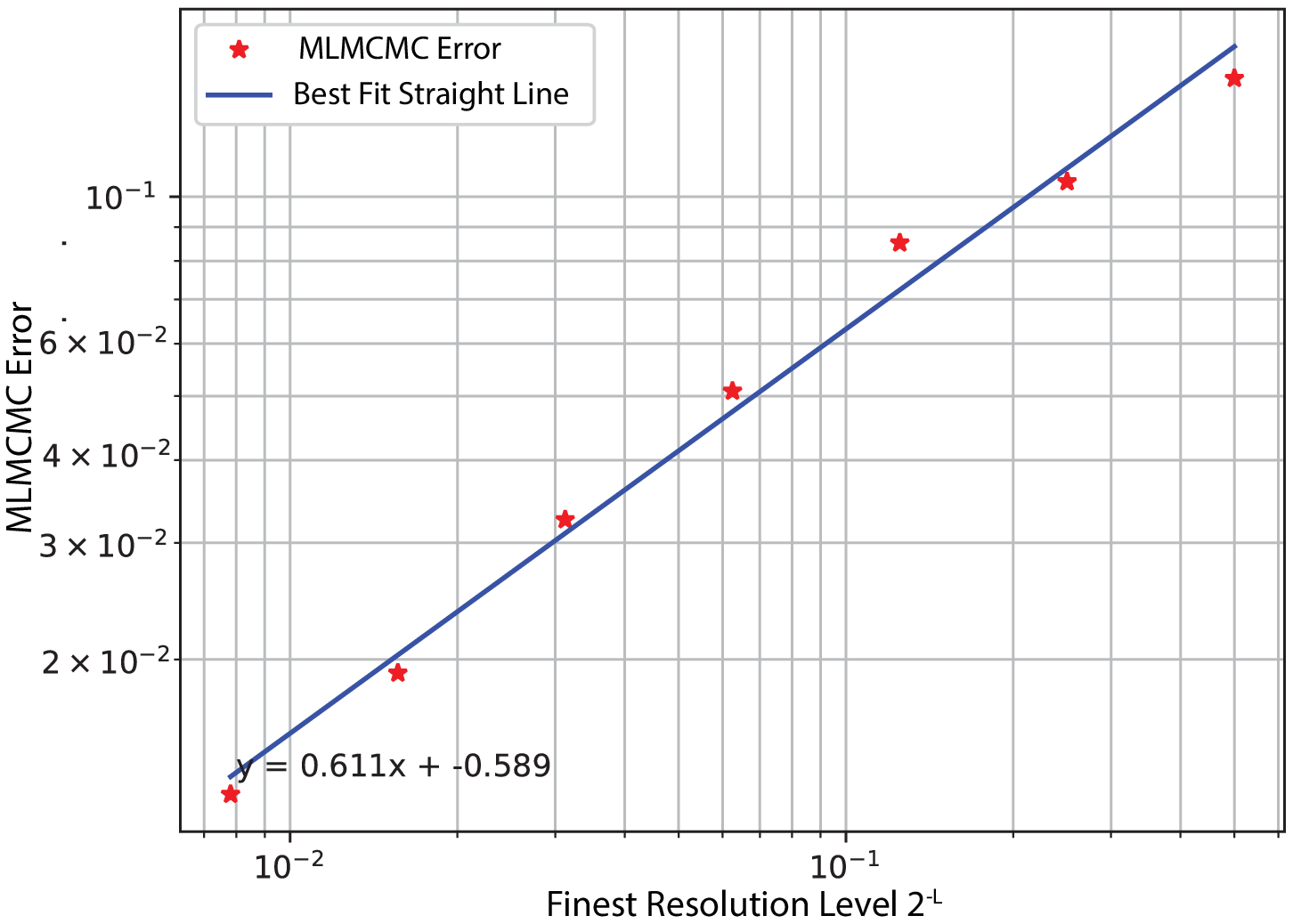}
    \caption{a = 4}
    \label{fig:pcn05alpha=4}
    \end{subfigure}
    \caption{pCN MLMCMC Errors for $\beta = 0.5$}
    \label{fig:pcn05mlmcmc}
\end{figure}

\begin{figure}[htbp!]
    \centering
    \begin{subfigure}{0.475\linewidth}
    \includegraphics[width=\linewidth]{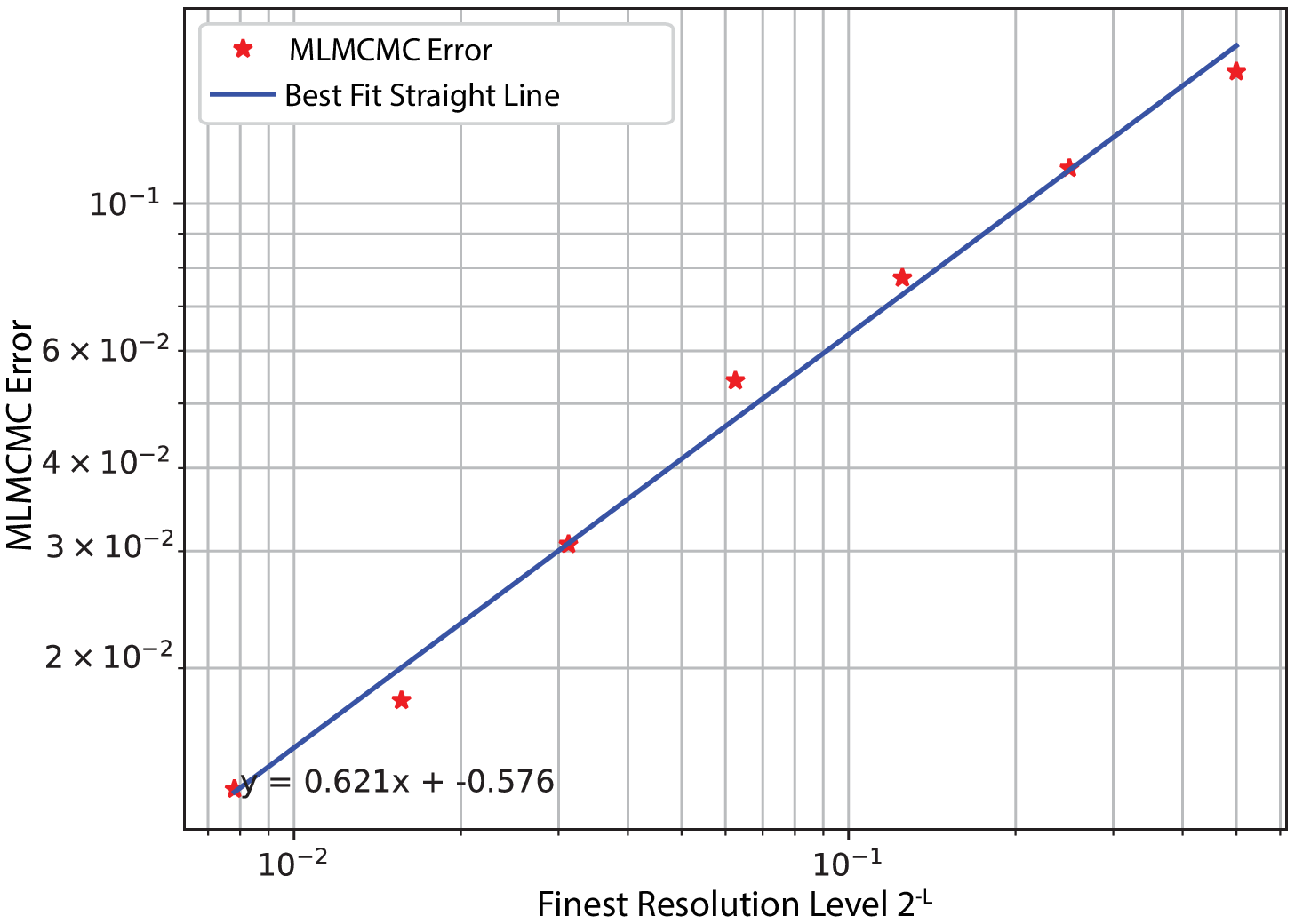}
   	\caption{a = 3}
    \label{fig:pcn02alpha=3}
    \end{subfigure}
    \begin{subfigure}{0.49\linewidth}
    \includegraphics[width=0.95\linewidth]{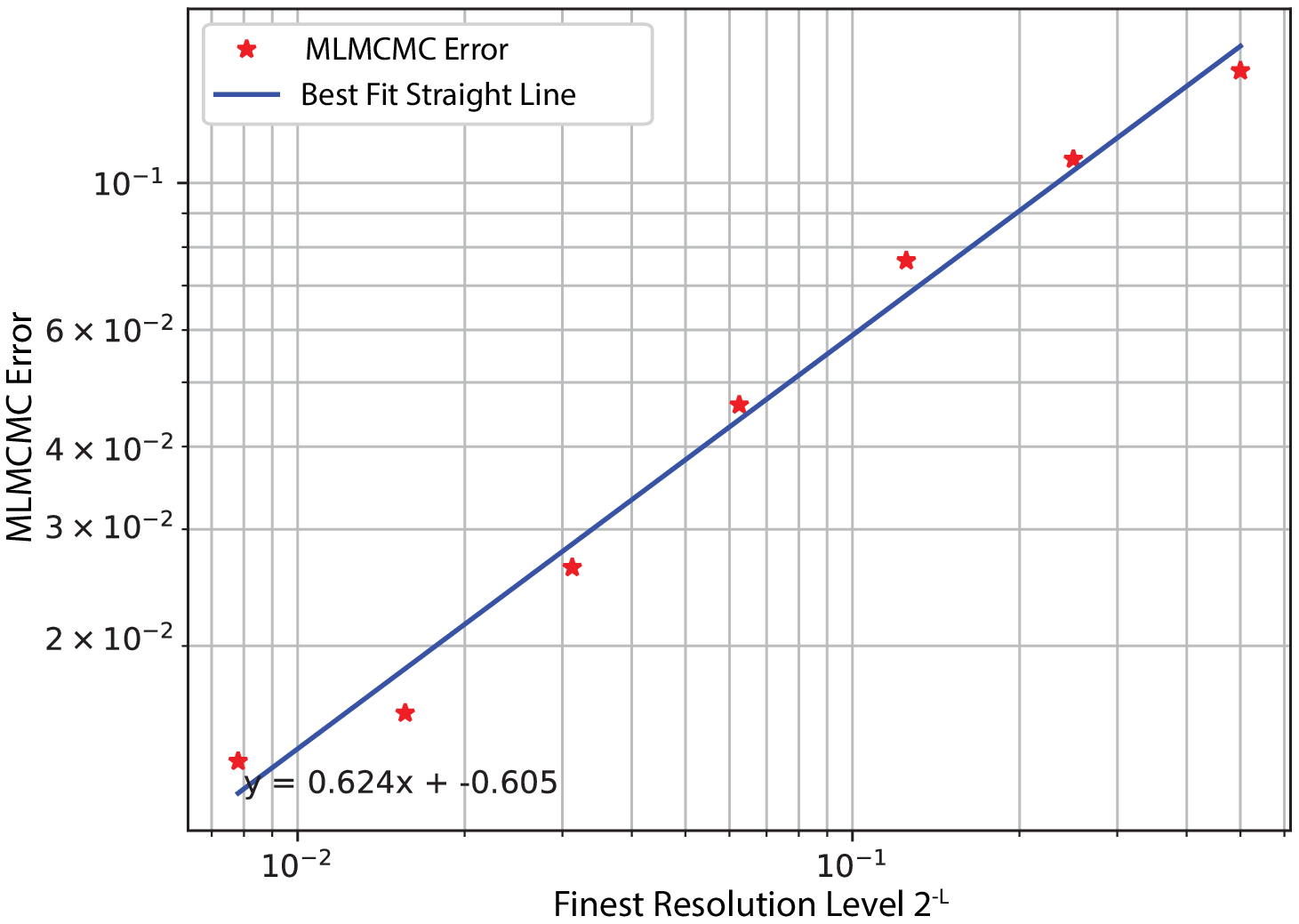}
    \caption{a = 4}
    \label{fig:pcn02alpha=4}
    \end{subfigure}
    \caption{pCN MLMCMC Errors for $\beta = 0.2$}
    \label{fig:pcn02mlmcmc}
\end{figure}

Figures 
\ref{fig:pcn08mlmcmc}, \ref{fig:pcn05mlmcmc}, and \ref{fig:pcn02mlmcmc} show the average absolute errors of the posterior expectation of the quantity of interest ($T((1/2, 1/2))$).
We use pCN sampler with the parameter value $\beta = 0.8, 0.5$, and $0.2$ respectively. The errors shown are absolute error averages of 32 independent MLMCMC runs. 
The best fit straight lines to the absolute errors of MLMCMC for $a = 3$ and $a = 4$ are 
0.538 and 0.581 for $\beta = 0.8$, 0.565 and 0.611 for $\beta = 0.5$, and 0.621 and 0.624 for $\beta = 0.2$, which are in reasonable agreement or slightly better than the theoretical convergence rate
in Table \ref{table:total error summary}.

\begin{figure}[htbp!]
    \centering
    \includegraphics[width=0.5\textwidth]{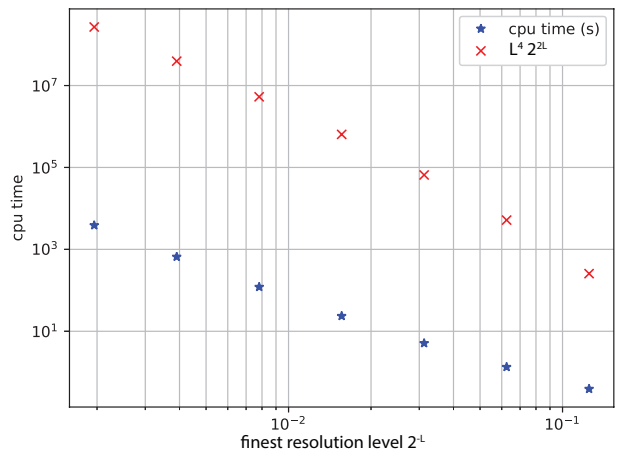}
 \caption{CPU time for a=3}
    \label{fig:cputimea=3}
\end{figure}

Next, in Figure \ref{fig:cputimea=3} we plot the  CPU time against the finest resolution grid size $2^{-L}$; the figure presents the average  CPU time of  five independent  runs of the MLMCMC algorithm, using independent sampler, where the heap sort procedure is employed in solving the forward eikonal equation by the FMM, for the case $a=3$ in \eqref{eq:a}. We plot also the bound $L^42^{2L}$ established in Remark \ref{rem:complexity} for comparison. The figure clearly indicates that the CPU time required is in agreement with the bound in Remark \ref{rem:complexity}.
%
%

For the case where the slowness depends on many random variables, we compare the MLMCMC results to the true solution profile of the forward eikonal equation with the reference slowness. Gauss-Hermite quadrature for computing an accurate reference posterior expectation is too complicated in this case. 
We  consider the multivariate slowness function of the form
\begin{equation}
\label{eq:logNslowness}
s({x}, u) = \exp(\sum_{i,j = 0}^{\infty} u_{ij} \phi_{ij} (x)),
\end{equation}
with $u_{i,j} \sim \NN (0, 1)$, 
 where
\begin{equation}
\label{eq:logNbasiscantor}
\phi_{i,j} ({x}) = \frac{\kappa}{((i+1)^2 + (j+1)^2)^2} \sin((i+1) \pi x_1) \sin((j+1) \pi x_2),
\end{equation}
for $x=(x_1,x_2)\in {\mathbb R}^2$;  $\kappa$ is a constant. 
We consider the eikonal equation in the two dimension domain $\Omega=(-1,1)\times(-1,1)$. This log-normal slowness can be recasted into the form of one index summation in \eqref{eq:lognormalslowness} where
\be
s(x,u)=\exp(\sum_{k=1}^\infty u_k\phi_k(x)),
\label{eq:sxunumerical}
\ee
with $u_k=u_{i_k,j_k}\sim \NN(0,1)$, and
\be
\phi_k(x)={\kappa\over ((i_k+1)^2+(j_k+1)^2)^2}\sin ((i_k+1)\pi x_1)\sin((j_k+1)\pi x_2),
\label{eq:phik}
\ee
for $x=(x_1,x_2)$, with the Cantor pairing function $k=\frac12(i_k+j_k)(i_k+j_k+1)+j_k+1$. Here the decaying rate in Assumption \ref{ass:1} $p=2$. 

\begin{figure}[htbp!]
\centering
\begin{subfigure}{0.32\linewidth}
\centering
\includegraphics[width=\linewidth]{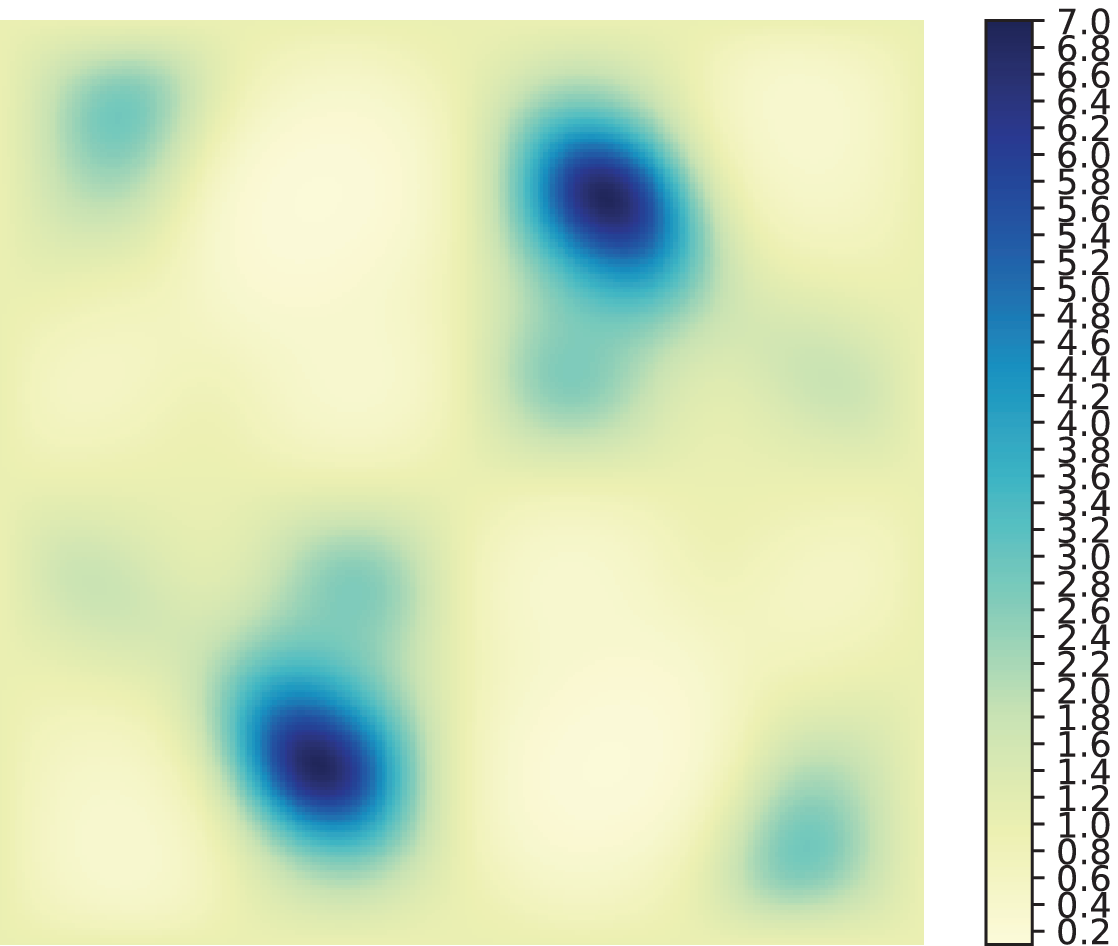}
\caption{Exact reference slowness}
\label{fig:symslw}
\end{subfigure}
\hspace{0.05cm}
\begin{subfigure}{0.32\linewidth}
\centering
\includegraphics[width=\linewidth]{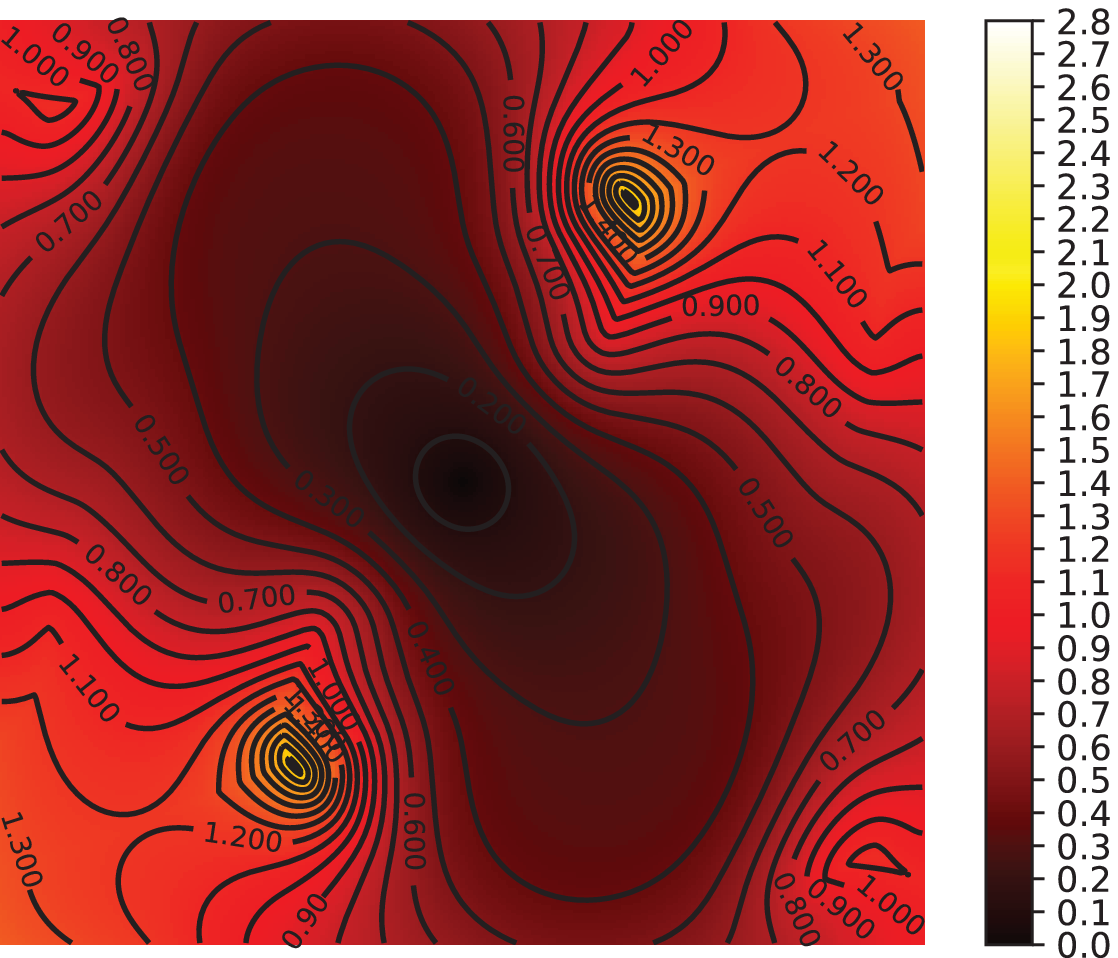}
\caption{Reference solution}
\label{fig:symTrefdouble}
\end{subfigure}
\hspace{0.05cm}
\begin{subfigure}{0.32\linewidth}
\centering
\includegraphics[width=\linewidth]{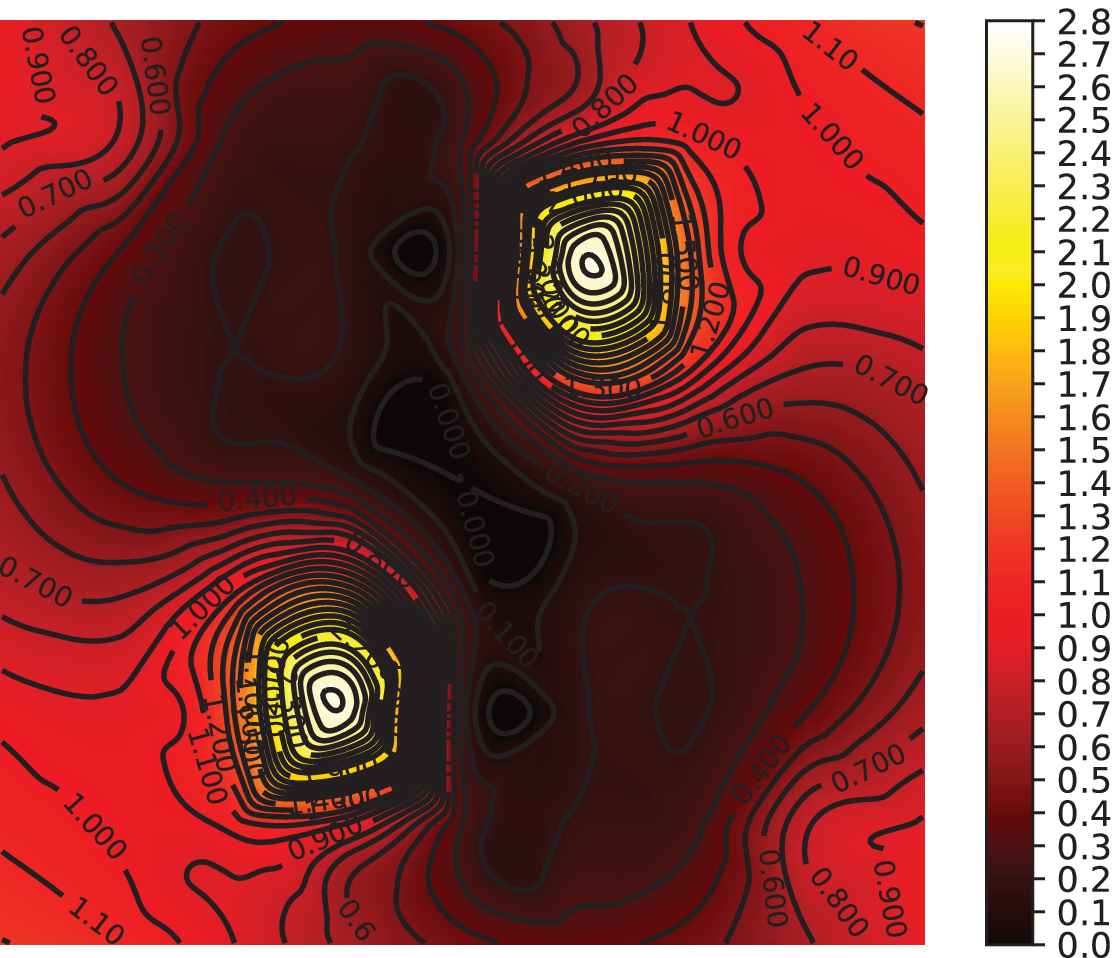}
\caption{MLMCMC recovered solution}
\label{fig:symTmlmcmcdouble}
\end{subfigure}
\hspace{0.05cm}
\caption{MLMCMC recovery of the solution}
\label{fig:symT}
\end{figure}

For the next experiment, to generate a reference slowness function, we take 64 terms in expansion \eqref{eq:logNslowness}, i.e., we choose $i,j=0,\ldots,7$. We generate a realization of the random log-normal slowness function as in Figure \ref{fig:symslw}, which contains two small areas with significantly larger slowness. In this experiment, we choose $\kappa=20$. 
For the Bayesian inverse problem, the observation is the solution of the forward eikonal equation at 64 points which are of equidistance of 1/8 on the boundary of $\Omega$. To generate the observation, we use data from the solution of the forward eikonal equation \eqref{eq:eikonal}, obtained from FMM with mesh size $O(2^{-12})$, at these 64 boundary points, for  5 different source points $x_0$ in $\Omega$ which are evenly spread throughout the domain $\Omega$. The observation vector is thus of dimension 320. The noisy observation $\delta$ in \eqref{eq:delta} is obtained by adding to this observation data a randomly generated realization of the normal distribution $\NN(0, 10^{-4}I)$ where $I$ is the $320\times 320$ dimension identity matrix. We use these data to recover the solution of the forward equation \eqref{eq:eikonal} for the source point $x_0=(0,0)$. The  true reference solution $T$ profile in Figure \ref{fig:symTrefdouble} is generated by solving the forward eikonal equation with the fine mesh $2^{-12}$. For the MLMCMC sampling procedure in Section \ref{sec:mlmcmc}, the coarsest level $l_0=4$ and the finest level $L=8$. The parameter $a$ in \eqref{eq:a} is 3. For the coarsest level of sampling where $l=l'=l_0$, where the forward equation is cheap to solve, we use 10000 samples, after discarding 5000 samples for burning in. For other discretizing levels, the number of samples are  as in Section \ref{sec:mlmcmc}, with the initial sample of the Markov chain is taken as the average of the samples of the previous level. The recovery solution of the forward eikonal equation is depicted in Figure \ref{fig:symTmlmcmcdouble}, which is the average of 8 independent runs of the MLMCMC. For each run of the MLMCMC, we compute the approximated posterior expectation of the solution of the forward eikonal equation at the nodes of a $8\times 8$ uniform grid in $\Omega$, and use bicubic interpolation to approximate the solution at other points. We note that Figure \ref{fig:symTrefdouble} depicts the true solution, while Figure \ref{fig:symTmlmcmcdouble} depicts the approximation to the posterior expectation of the solution obtained by MLMCMC. In general, the true solution for a reference slowness is not equal to the posterior expectation of the forward solution. However, we see clearly that the posterior expectation of the solution obtained from MLMCMC is in a reasonable agreement with the true solution. It detects accurately the two areas of maximum value of the slowness, and also the two areas of a higher value of the slowness in the top left and bottom right corners. 

In Section \ref{sec:mlmcmc}, the quantity of interest is the forward solution. 
The MLMCMC algorithm is capable of approximating the posterior expectation of the slowness in a similar fashion. We consider the prior log-normal form of the slowness in \eqref{eq:logNslowness} and \eqref{eq:logNbasiscantor}.  The function $\phi_k$ in  \eqref{eq:phik} has the decay rate $p=2$ in Assumption \ref{ass:1}. Choosing the truncation level $J_{l'}=O( 2^{l'/2})$ in \eqref{eq:truncatedslowness}, we denote by
\[
s^{l'}(x,u)=\exp(\sum_{k=1}^{J_{l'}}u_k\phi_k(x)).
\]
We have
\[
|s(x,u)-s^{l'}(x,u)|\le c \left(\sum_{i=1}^{J_{l'}}|u_i|b_i\right)\exp(\sum_{i=1}^\infty |u_i|b_i).
\]
The MLMCMC algorithm for approximating the posterior expectation of the slowness, with the log-normal prior, is obtained by choosing the function $Q(u)$ in Section \ref{sec:mlmcmc} by $Q(u)=s^{l'}(x,u)-s^{l'-1}(x,u)$ for $l'\ge l_0+1$, and $Q(u)=s^{l_0}(x,u)$ for $l'=l_0$. We note that this bears some resemblance to the QMC approximation of the solution of the forward random elliptic equation in \cite{Dick2019}. 

\begin{figure}[htbp!]
\centering
\begin{subfigure}{0.44\linewidth}
\centering
\includegraphics[width=\linewidth]{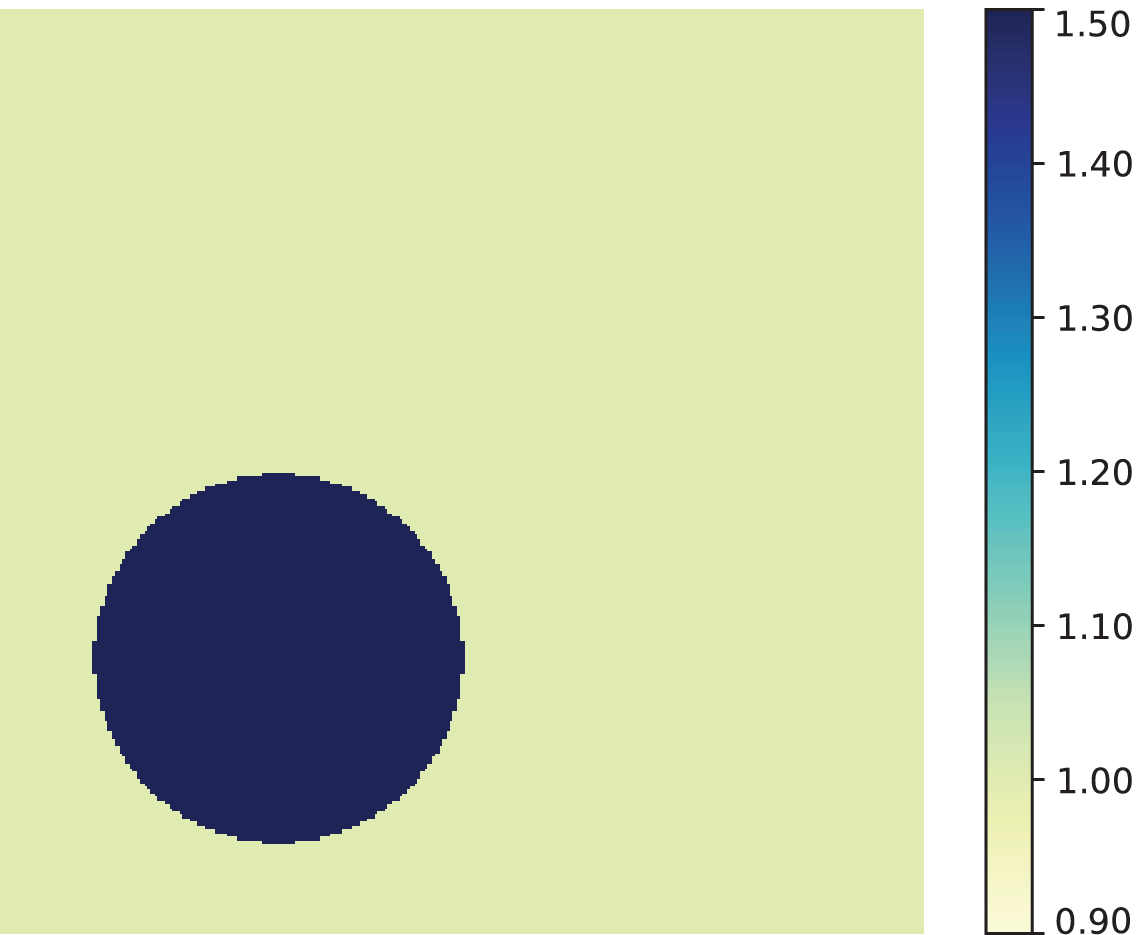}
\caption{Reference slowness}
\label{fig:kappatuning}
\end{subfigure}
\hspace{0.1cm}
\begin{subfigure}{0.45\linewidth}
\centering
\includegraphics[width=\linewidth]{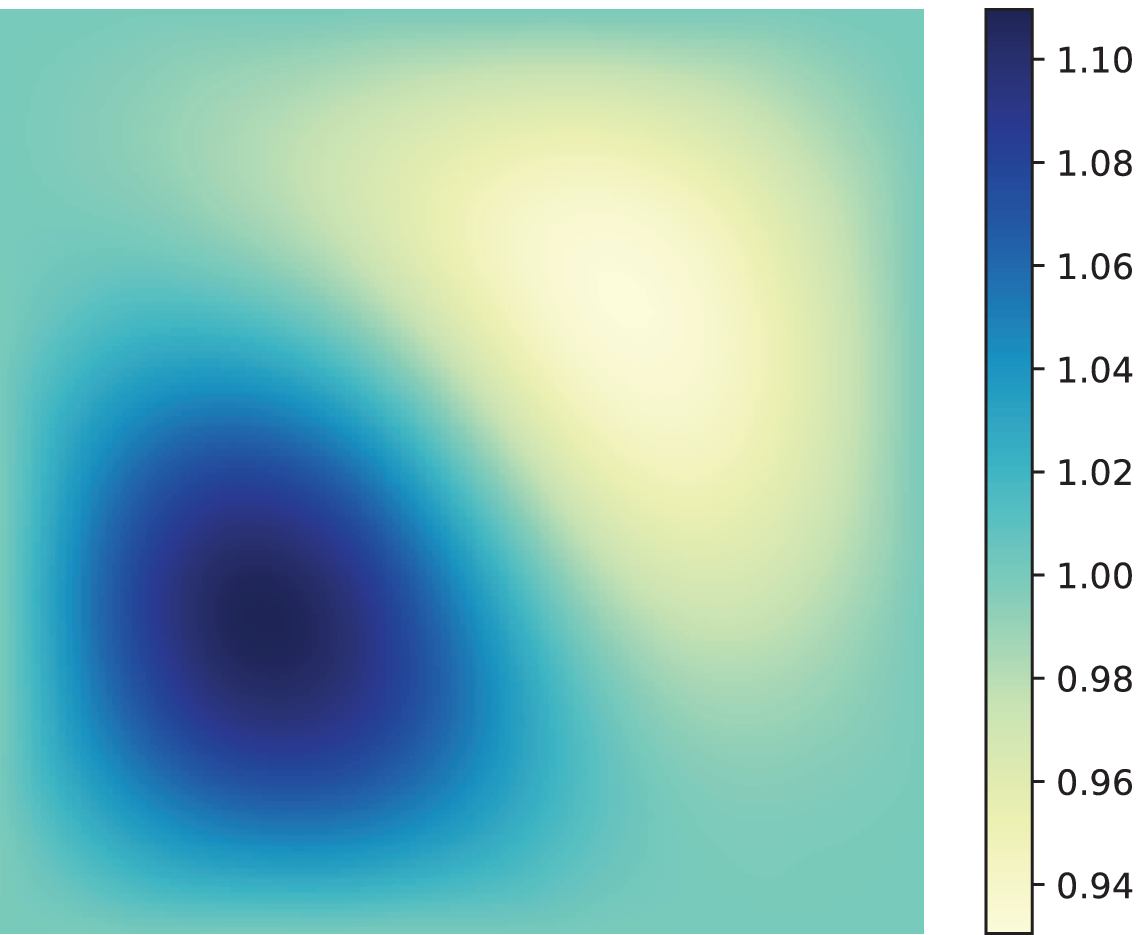}
\caption{MLMCMC recovery of the slowness for $\kappa = 1$}
\label{fig:kappa1tune}
\end{subfigure}
\hspace{0.1cm}
\begin{subfigure}{0.45\linewidth}
\centering
\includegraphics[width=\linewidth]{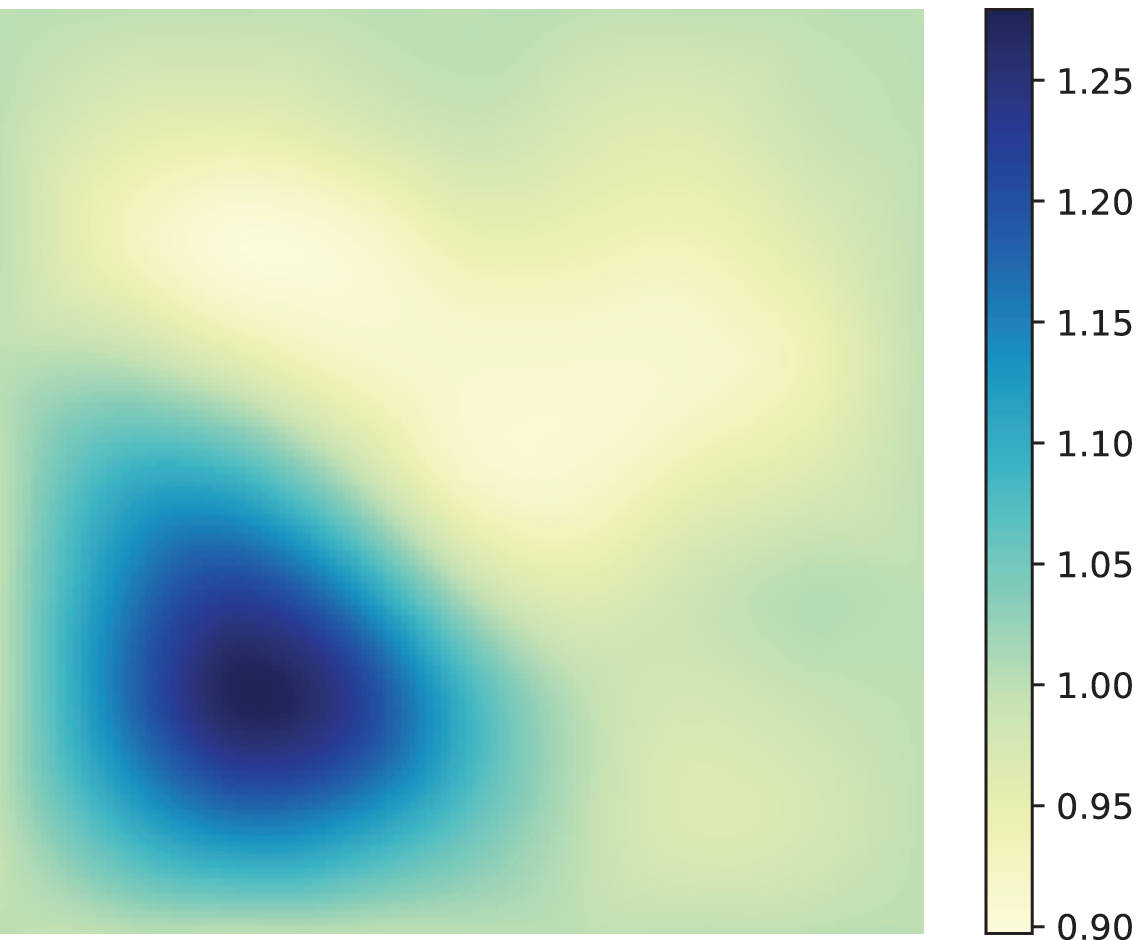}
\caption{MLMCMC recovery of the slowness for $\kappa = 10$}
\label{fig:kappa10tune}
\end{subfigure}
\hspace{0.1cm}
\begin{subfigure}{0.45\linewidth}
\centering
\includegraphics[width=\linewidth]{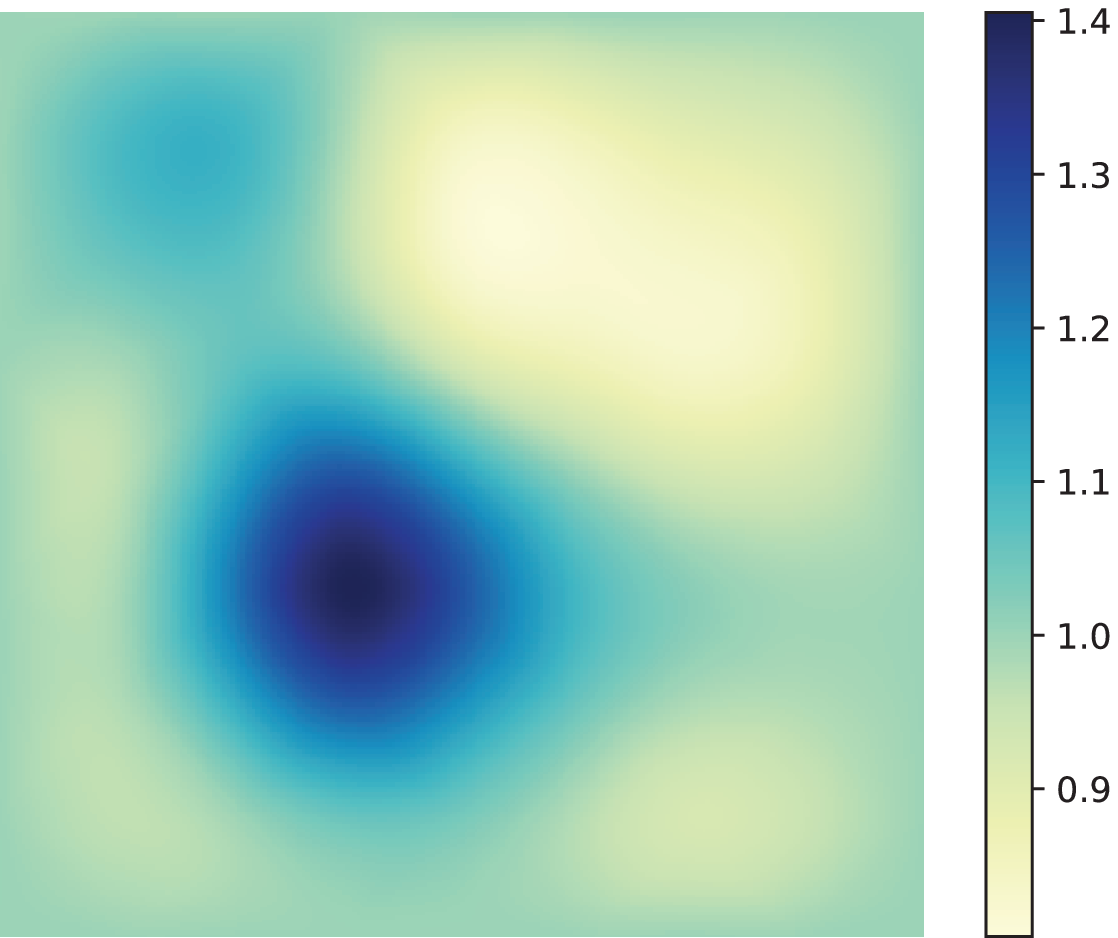}
\caption{MLMCMC recovery of the slowness for $\kappa = 20$}
\label{fig:kappa20tune}
\end{subfigure}
\hspace{0.1cm}
\caption{MLMCMC recovery of the binary slowness; slowness value 1.5 and 1}
\label{fig:tuning}
\end{figure}

\begin{figure}[htbp!]
\centering
\begin{subfigure}{0.32\linewidth}
\centering
\includegraphics[width=\linewidth]{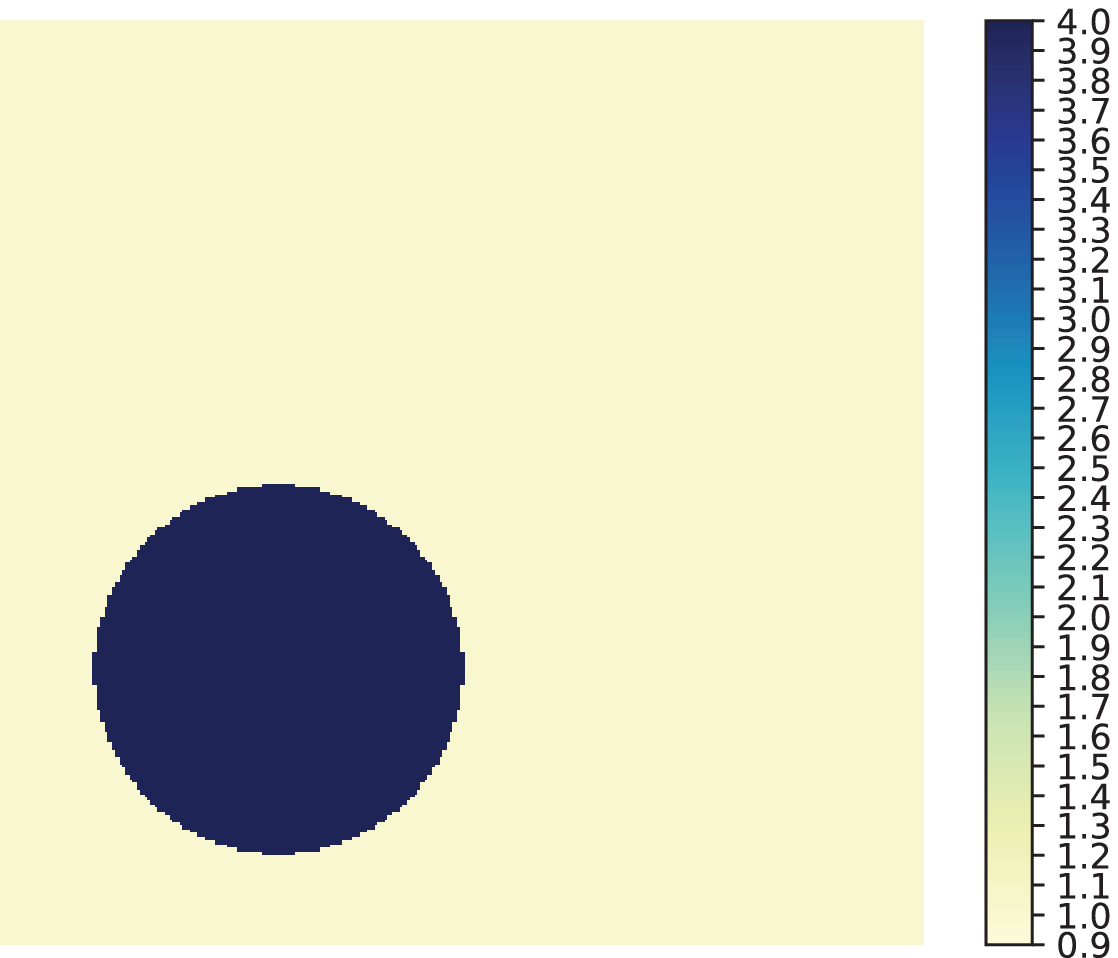}
\caption{Reference slowness}
\label{fig:k1denseorig}
\end{subfigure}
\hspace{0.05cm}
\begin{subfigure}{0.32\linewidth}
\centering
\includegraphics[width=\linewidth]{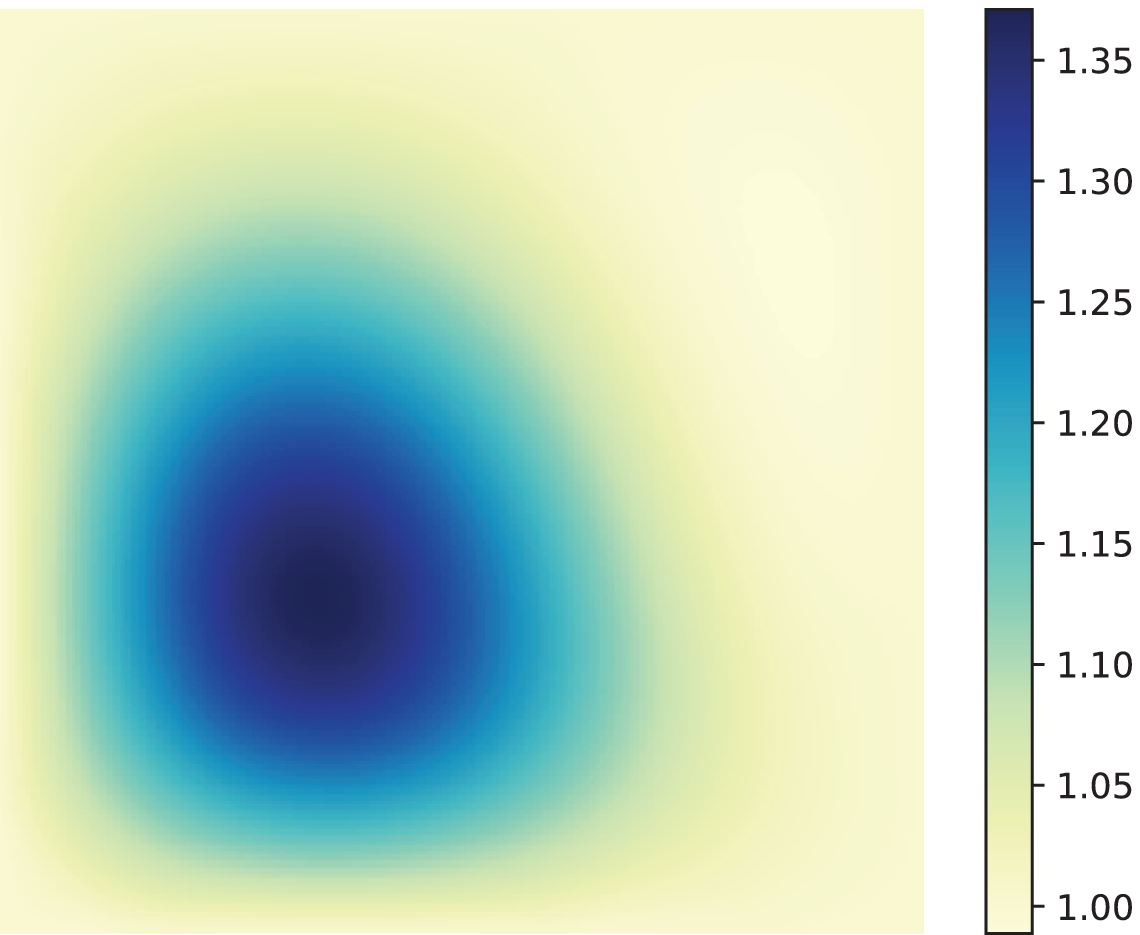}
\caption{MLMCMC recovery of the slowness for $\kappa = 1$}
\label{fig:k1densedoubleraw}
\end{subfigure}
\hspace{0.05cm}
\begin{subfigure}{0.32\linewidth}
\centering
\includegraphics[width=\linewidth]{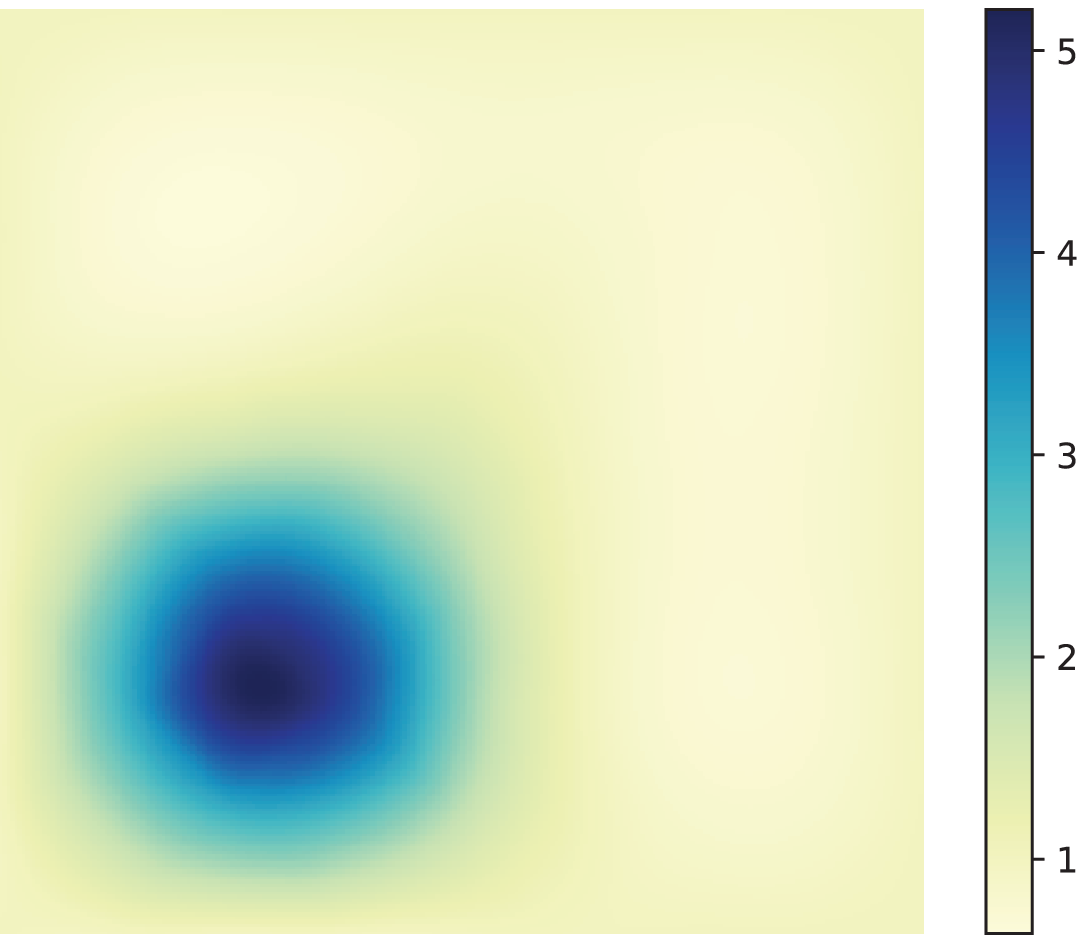}
\caption{MLMCMC recovery of the slowness for $\kappa = 20$}
\label{fig:k20densedoubleraw}
\end{subfigure}
\hspace{0.05cm}
\caption{MLMCMC recovery of the binary slowness; slowness value 4 and 1}
\label{fig:porousinclusion}
\end{figure}

\begin{figure}[htbp!]
\centering
\begin{subfigure}{0.32\linewidth}
\centering
\includegraphics[width=\linewidth]{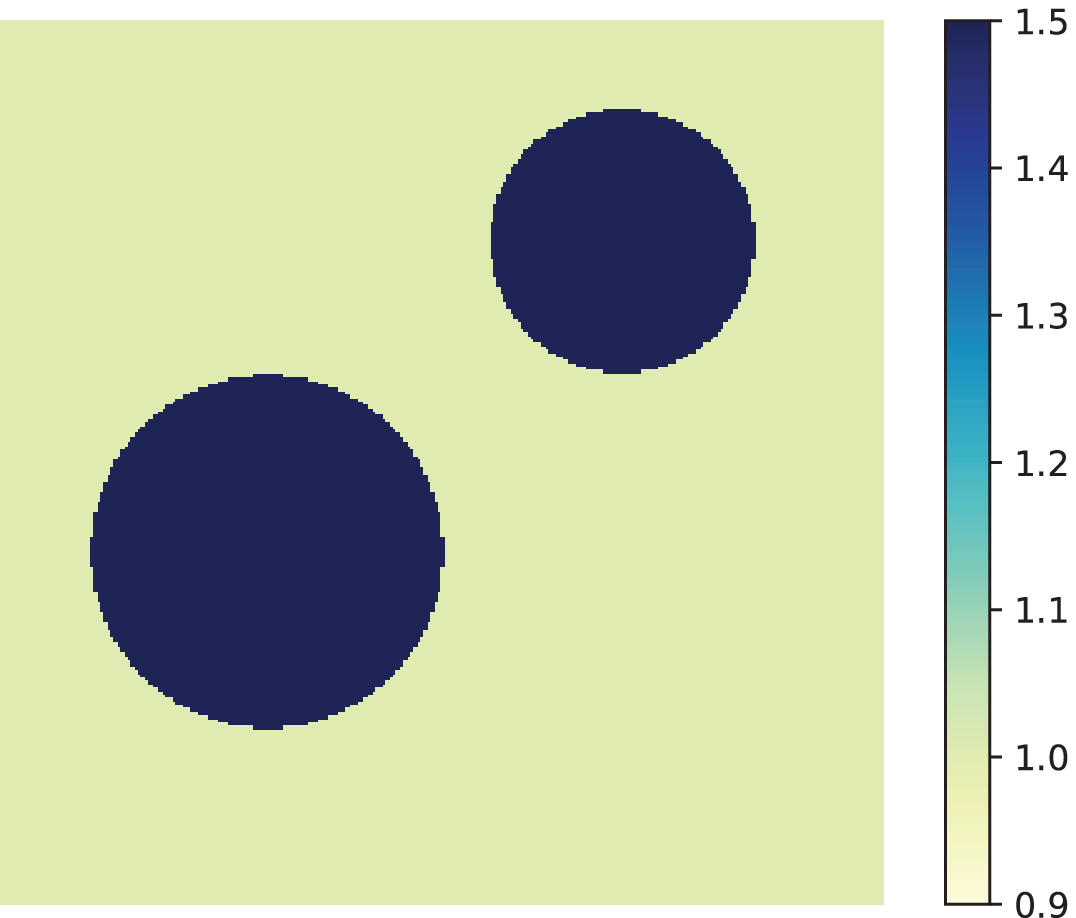}
\caption{Reference slowness}
\label{fig:k1dualorig}
\end{subfigure}
\hspace{0.05cm}
\begin{subfigure}{0.32\linewidth}
\centering
\includegraphics[width=\linewidth]{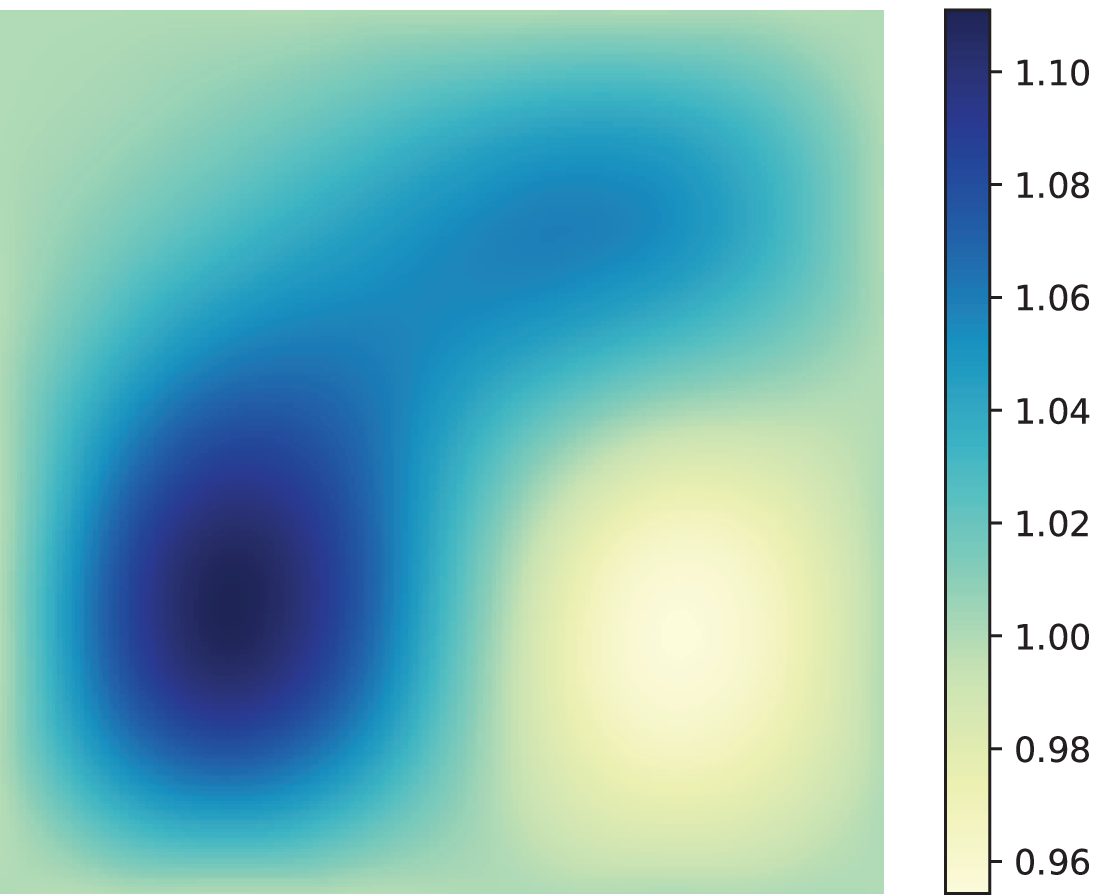}
\caption{MLMCMC recovery of the slowness for $\kappa = 1$}
\label{fig:k1dualdoubleraw}
\end{subfigure}
\hspace{0.05cm}
\begin{subfigure}{0.32\linewidth}
\centering
\includegraphics[width=\linewidth]{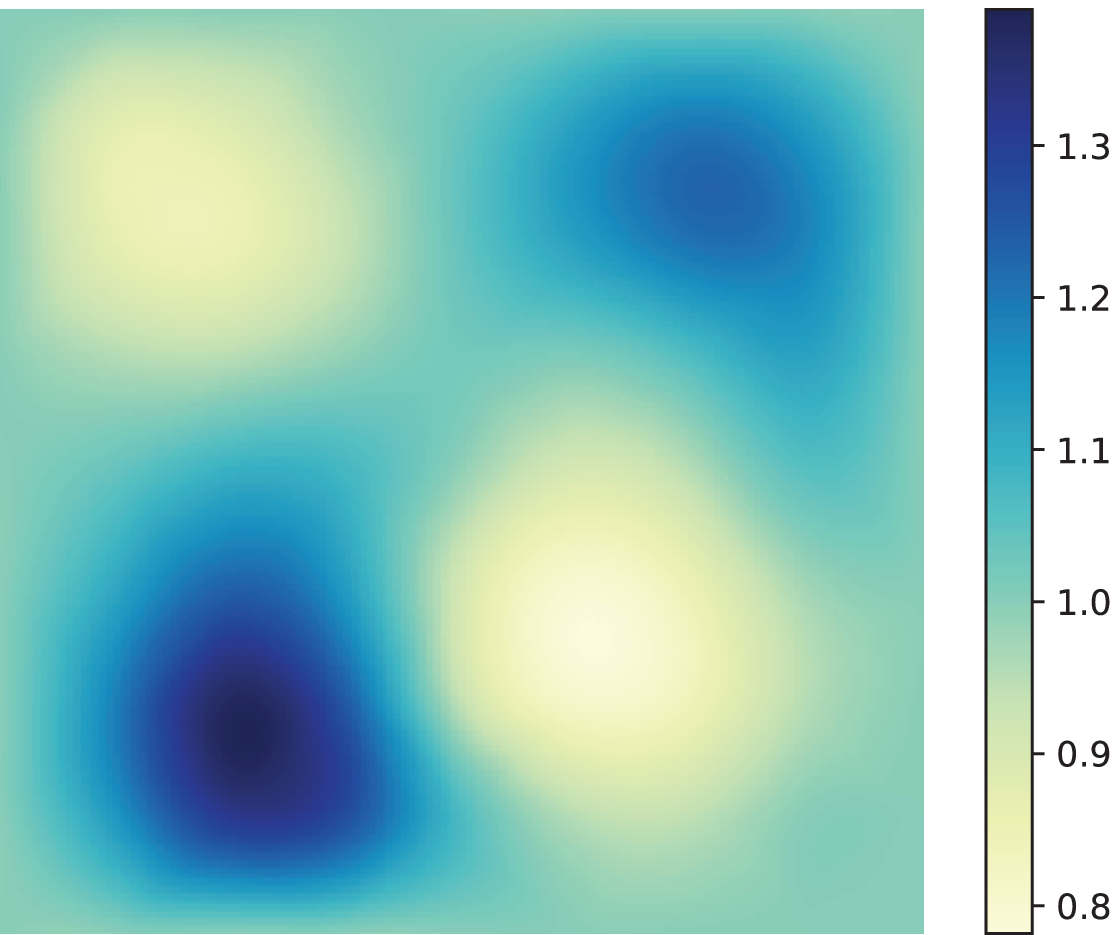}
\caption{MLMCMC recovery of the slowness for $\kappa = 20$}
\label{fig:k20dualdoubleraw}
\end{subfigure}
\hspace{0.05cm}
\caption{MLMCMC recovery of two inclusions}
\label{fig:dualinclusion}
\end{figure}

In the previous numerical examples, we generate a reference slowness function by taking a realization of the prior log-normal slowness. Now we generate observation data from a slowness function which is not a priori related to the log-normal form. We consider the binary slowness in Figure \ref{fig:kappatuning}. The domain $\Omega$ in this case is $(0,1)\times (0,1)$.  For the next experiment, we choose the value of the slowness inside the circle inclusion to be 1.5, and the value of the slowness outside to be 1. The data are generated similarly to the last example. We  take the solution of the forward eikonal equation at 64 equi-distanced boundary points. The observation data is generated by solving the forward eikonal equation \eqref{eq:eikonal} 5 times for 5 different source points $x_0$, with small mesh size $2^{-12}$, making it a 320 dimensional vector. To generate the noisy observation $\delta$ in \eqref{eq:delta}, we add a randomly generated realization of the noise distribution $N(0,10^{-4}I)$ where $I$ is the $320\times 320$  dimensional identity matrix. We choose $J_l=2\lceil 2^{l/2}\rceil$ and $J_{l'}=2\lceil 2^{l'/2}\rceil$.
 Figure \ref{fig:tuning} presents the recovery of the slowness for $\kappa=1,10$ and $20$, where the results is the average of the outputs of 8 independent runs of the MLMCMC algorithm. The results demonstrate that the MLMCMC sampling procedure in Section \ref{sec:mlmcmc} for the Bayesian inverse problem with the log-normal prior in \eqref{eq:logNslowness} and \eqref{eq:logNbasiscantor}, using data from the binary slowness, can recover fairly accurately the position of the inclusion of higher slowness. However, the value of the slowness inside the inclusion is captured more accurately with a higher value of $\kappa$. This is because a typical realization of the log-normal slowness has a larger absolute value when $\kappa$ is larger. We observe the same outcome in Figure \ref{fig:porousinclusion}, where the value of the slowness inside the inclusion is now changed to 4. We can always recover the area of the inclusion, but we recover accurately the high value of the slowess inside the inclusion when $\kappa$ is larger. Figure \ref{fig:dualinclusion} presents the recovery results for the case of two inclusions. The figure again shows that MLMCMC is capable of detecting the inclusions. 

\begin{figure}[htbp!]
\centering
\begin{subfigure}{0.32\linewidth}
\centering
\includegraphics[width=\linewidth]{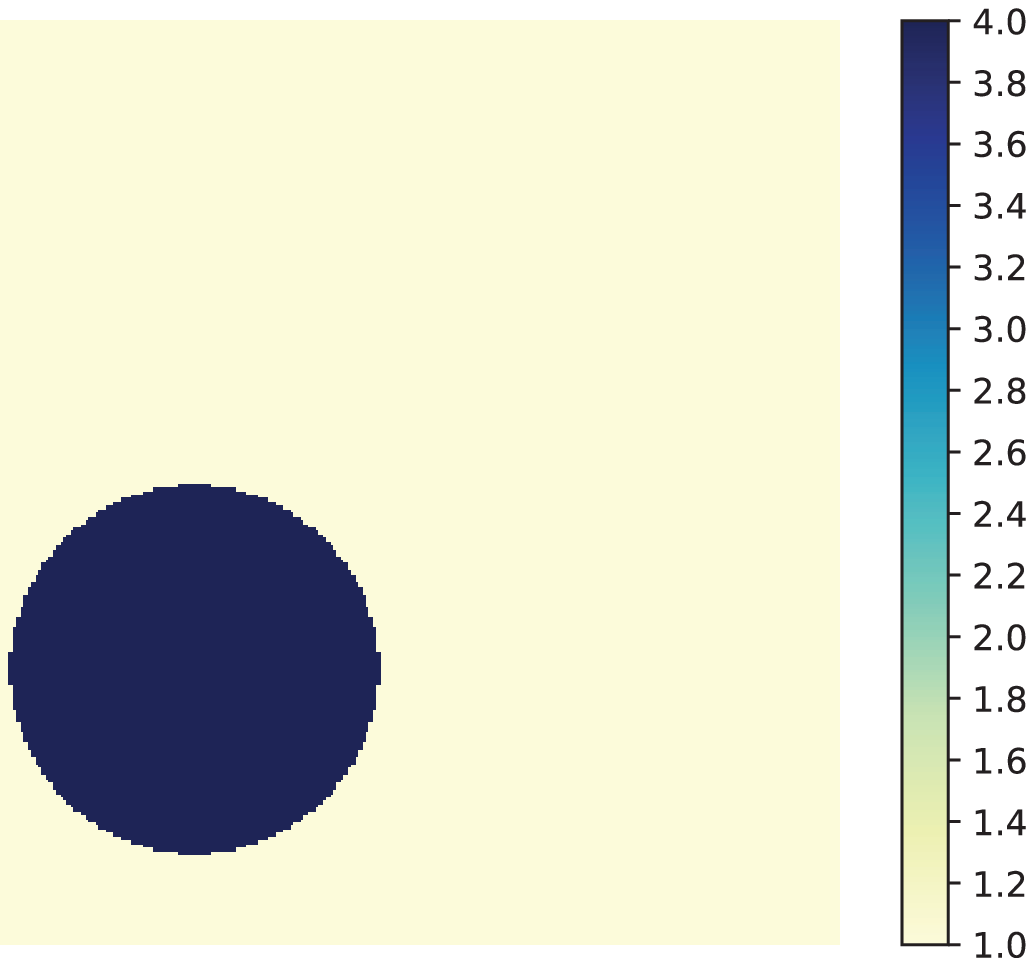}
\caption{Exact reference slowness}
\label{fig:binTdenseslw}
\end{subfigure}
\hspace{0.05cm}
\begin{subfigure}{0.32\linewidth}
\centering
\includegraphics[width=\linewidth]{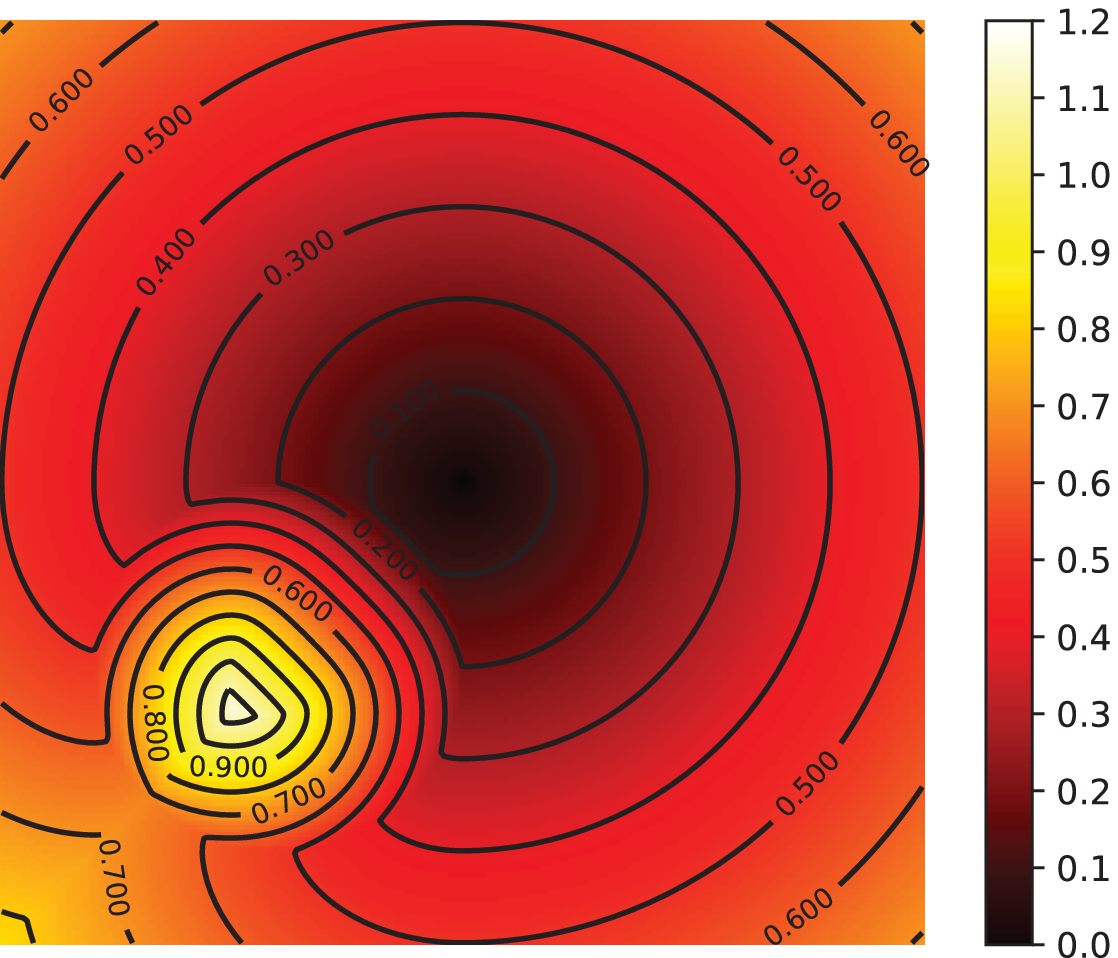}
\caption{Reference solution $T$}
\label{fig:binTdenserefdouble}
\end{subfigure}
\hspace{0.05cm}
\begin{subfigure}{0.32\linewidth}
\centering
\includegraphics[width=\linewidth]{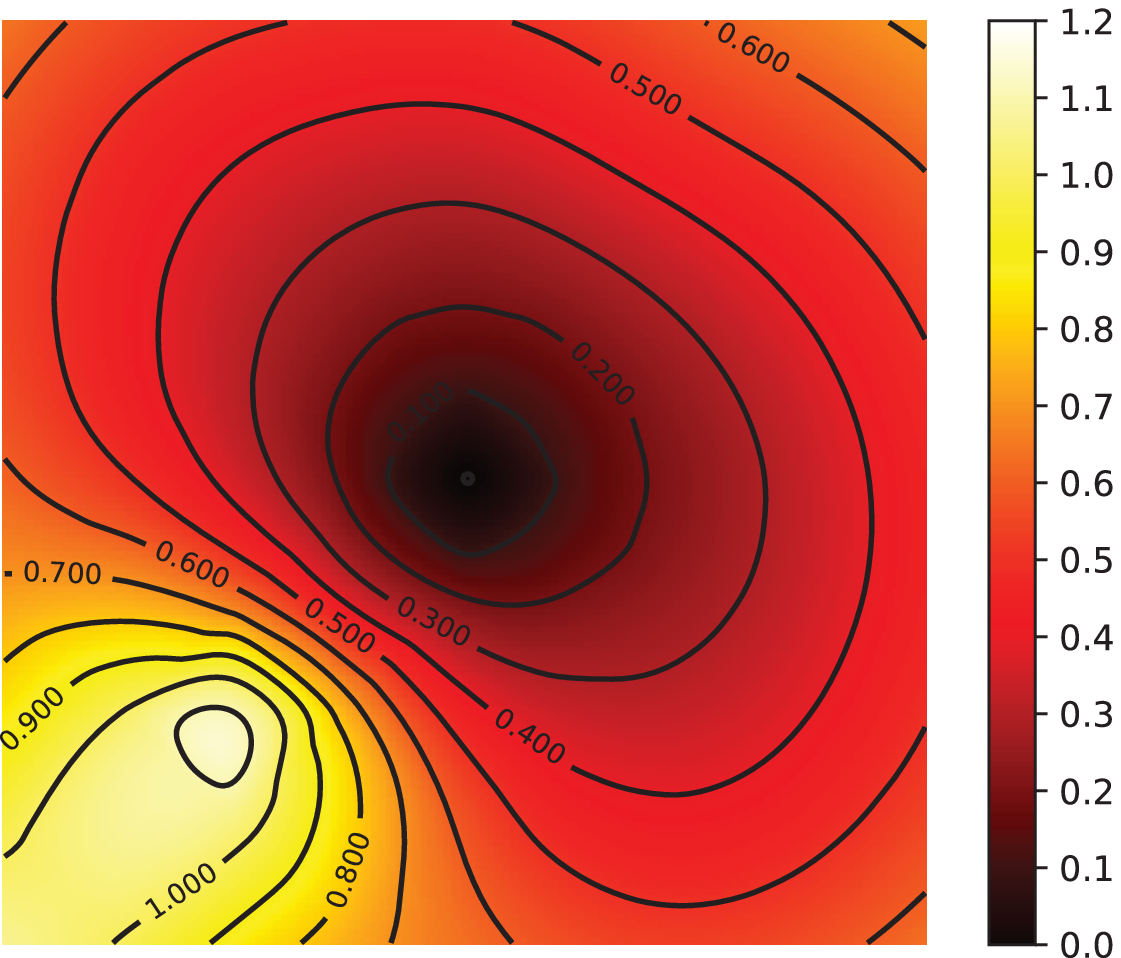}
\caption{MLMCMC recovered solution}
\label{fig:binTdensemlmcmcdouble}
\end{subfigure}
\hspace{0.05cm}
\caption{MLMCMC recovery for the forward solution for binary slowness; slowness value 4 and 1}
\label{fig:binTdense}
\end{figure}
In Figure \ref{fig:binTdense}, we present the recovery of the solution to the forward equation where the data are generated from  the reference binary slowness in Figure \ref{fig:porousinclusion} with the value 4 of the slowness inside the circle inclusion, and value  1 outside. The  log-gaussian prior is in \eqref{eq:logNslowness} and \eqref{eq:logNbasiscantor}.
We use MLMCMC to compute the posterior expectation of the solution at the nodes of a uniform $8\times 8$ grid inside the domain $\Omega=(0,1)\times (0,1)$. We use bicubic interpolation to approximate the posterior expectation of the forward solution at other points in $\Omega$. 
We depict the true solution for the forward equation, and the average results of 8 independent runs of the MLMCMC.
The MLMCMC algorithm is clearly able to  produce the  shock front in the forward solution, which is due to the inclusion of significantly higher slowness value. 

\vskip 30pt
{\bf Acknowledgement} Zhan Fei Yeo's research is supported by the Nanyang President's Graduate Scholarship. Viet Ha Hoang's research is supported by the Singapore
Ministry of Education Academic Research Fund Tier 2 grant MOE2017-T2-2-144. We thank Dr Jia Hao Quek for assisting with the implementation of MLMCMC at the beginning of the project.

\bibliographystyle{plain}
\bibliography{citation}

\appendix

\section{Some inequalities}\label{app:A}
%
\begin{lem}
	\label{lem:exponentialbounds}Let $t \in \R^+$. We have the following inequalities.
	
\begin{equation}
		\label{eq:exp1}
		\int_{-\infty}^{\infty} \exp( - \frac{z^2}{2} + \abs{z}t) \frac{dz}{\sqrt{2\pi}} \leq \exp(\frac{ t^2}{2})\exp(t \sqrt{\frac{2}{\pi}}),
		\end{equation}
\begin{equation}\label{eq:exp2}
		\int_{-\infty}^{\infty} z^2\exp(- \frac{z^2}{2} + \abs{z}t) \frac{dz}{\sqrt{2\pi}} \leq C\exp(\frac{t^2}{2})(1 + t^2),
		\end{equation}
\begin{equation}\label{eq:exp3}
		\int_{-\infty}^{\infty} \abs{z} \exp(- \frac{z^2}{2} + \abs{z}t) \frac{dz}{\sqrt{2\pi}} \leq C\exp(\frac{t^2}{2})\left( 1+t\right).
		\end{equation}
\end{lem}
We refer to \cite{HoangSchwabmcmclogn} for the proofs of these inequalities. 

\section{Justification of MLMCMC convergence rates}
\label{app:justification}
It is necessary to modify the rigorous proof in \cite{hoang2020analysis} of the convergence rate of the MLMCMC method in Section \ref{sec:mlmcmc}  as the theoretical convergence rate of the FMM for the eikonal equation in \eqref{eq:fmmrate} only holds when the grid size $h$ is not more than an upper bound $h_0(u)$, which can be arbitrarily small for different realizations $u$. We present the main modifications in this appendix. Also the theoretical convergence rate of the FMM method is only $O(h^{1/2})$, which is weaker than the $O(h)$ convergence rate for the finite element method in \cite{hoang2020analysis}.

We note that 
\[
|1-\exp(\Phi^l(u;\delta)-\Phi^{l-1}(u;\delta))|\le |\Phi^l(u;\delta)-\Phi^{l-1}(u;\delta)||1+\exp(\Phi^l(u;\delta)-\Phi^{l-1}(u;\delta))|.
\]
Thus
\[
|A_1^{ll'}|\le c|(\Phi^l(u;\delta)-\Phi^{l-1}(u;\delta))(T^{l'}(x^*,u)-T^{l'-1}(x^*,u))|.
\]
We have further that
\beqas
|\Phi^l(u; \delta)) -\Phi^{l-1} (u; \delta))| &\leq c |2\delta - G^l(u) - G^{l-1} (u)|_\Sigma |G^l(u) - G^{l-1}(u)|_\Sigma\\
	&\leq c (|\delta|_\Sigma + |G^l(u)|_\Sigma + |G^{l-1} (u)|_\Sigma)|G^l (u) - G^{l-1}(u)|_\Sigma.
\eeqas
Let $U_1^{ll'}\subset U$ be the set of $u\in U$ such that $2^{-(l-1)}\le h_0(u)$ and $2^{-(l'-1)}\le h_0(u)$, i.e. the FMM convergence rate \eqref{eq:fmmrate} holds for the grid sizes $2^{-(l-1)}, 2^{-l}, 2^{-(l'-1)}$ and $2^{-l'}$. For $u\in U_1^{ll'}$, from \eqref{eq:boundGJhu},  \eqref{eq:truncationrate} and \eqref{eq:fmmrate}, we have
\be
|T^{l'}(x^*,u)-T^{l'-1}(x^*,u)|\le c\exp(c\sum_{i=1}^\infty(b_i+\bar b_i)|u_i|) (2^{-l'}+\sum_{i>J_{l'-1}}|u_i|b_i),
\label{eq:1}
\ee
and 
\be
|G^l(u)-G^{l-1}(u)| \le c\exp(c\sum_{i=1}^\infty(b_i+\bar b_i)|u_i|) (2^{-l}+\sum_{i>J_{l-1}}|u_i|b_i).
\label{eq:2}
\ee
We thus have
\[
|A_1^{ll'}|\le c\exp(c\sum_{i=1}^\infty |u_i|(b_i+\bar b_i))(2^{-l}+\sum_{i> J_{l-1}}|u_i|b_i)(2^{-l'}+\sum_{i>J_{l'-1}}|u_i|b_i).
\]
for $u\in U_1^{ll'}$. 
Let $\varepsilon=\sum_{i>J_{l-1}}b_i$ and $\varepsilon'=\sum_{i>J_{l'-1}}b_i'$. Using the inequality $x\le \varepsilon\exp(\frac{x}{\varepsilon})$, we have
\[
|A_1^{ll'}|\le c2^{-(l+l')/2}\exp(c\sum_{i=1}^\infty|u_i|(b_i+\bar b_i)+{1\over\varepsilon}\sum_{i>J_{l-1}}|u_i|b_i|+{1\over\varepsilon'}\sum_{i>J_{l'-1}}|u_i|b_i).
\]
Let $U_2^{ll'}\subset U$ be the set of $u\in U$ such that either $2^{-(l-1)\}}\le h_0(u)$ and $2^{-(l'-1)}> h_0(u)$ or $2^{-(l'-1)\}}\le h_0(u)$ and $2^{-(l-1)}> h_0(u)$, i.e. only one of the inequalities \eqref{eq:1} and \eqref{eq:2} holds. In this case, we have
\[
|A_1^{ll'}|\le c\exp(c\sum_{i=1}^\infty |u_i|(b_i+\bar b_i))(2^{-l}+\sum_{i> J_{l-1}}|u_i|b_i),
\]
or
\[
|A_1^{ll'}|\le c\exp(c\sum_{i=1}^\infty |u_i|(b_i+\bar b_i))(2^{-l'}+\sum_{i>J_{l'-1}}|u_i|b_i).
\]
Thus when $u\in U_2^{ll'}$, 
\[
|A_1^{ll'}|\le c2^{-\max\{l,l'\}/2}\exp(c\sum_{i=1}^\infty|u_i|(b_i+\bar b_i)+{1\over \min\{\varepsilon,\varepsilon'\}}\sum_{i>\max\{J_{l-1},J_{l'-1}\}}|u_i|b_i).
\]
Let $U_3^{ll'}=U\setminus (U_1^{ll'}\cup U_2^{ll'})$ be the set of $u\in U$ such that $2^{-\max\{l,l'\}}> h_0(u)$, i.e. neither \eqref{eq:1} nor \eqref{eq:2} hold. 
When $u\in U_3^{ll'}$, from \eqref{eq:boundTJh} and \eqref{eq:boundGJhu}, we have
\[
|A_1^{ll'}|\le c\exp(c\sum_{i=1}^\infty|u_i|b_i).
\]
Thus for all $u\in U$,
\begin{eqnarray}
&&|A_1^{ll'}|\le  c2^{-(l+l')/2}\exp(c\sum_{i=1}^\infty|u_i|(b_i+\bar b_i)+{1\over\varepsilon}\sum_{i>J_{l-1}}|u_i|b_i|+{1\over\varepsilon'}\sum_{i>J_{l'-1}}|u_i|b_i)+\\ \nonumber
&&\qquad c2^{-\max\{l,l'\}/2}\exp(c\sum_{i=1}^\infty|u_i|(b_i+\bar b_i)+{1\over \min\{\varepsilon,\varepsilon'\}}\sum_{i>\max\{J_{l-1},J_{l'-1}\}}|u_i|b_i){\cal I}_{U_2^{ll'}}(u)+ c\exp(c\sum_{i=1}^\infty|u_i|b_i){\cal I}_{U_3^{ll'}}(u),
\label{eq:Vll}
\end{eqnarray}
where ${\cal I}$ is the indicator function. We have similar estimates for other terms $A_j^{ll'}$ for $j=2,4,6,7,8$ and $A_j^l$ for $j=3,5$. First we consider the independence sampler. Let ${\cal E}^{\bar\gamma}$ be the expectation over the space of all the Markov chains generated by the MCMC process with the acceptance probability \eqref{eq:acceptanceprob}, and with the initial sample $u^{(0)}$ being distributed according to the probability $\bar\gamma$. We then have the following result.
\begin{lem}
For $g\in L^2(U,\gamma)$, 
\[
{\cal E}^{\bar\gamma}\left[\Big|{1\over M}\sum_{k=1}^Mg(u^{(k)})-\mathbb{E}^{\gamma^l}[g]\Big|^2\right]\le cM^{-1}\mathbb{E}^{\gamma}[g^2].
\]
\end{lem}
We refer to \cite{hoang2020analysis} Lemma B2 for a proof. We note that
\[
\int_U\left[\exp(c\sum_{i=1}^\infty|u_i|(b_i+\bar b_i)+{1\over\varepsilon}\sum_{i>J_{l-1}}|u_i|b_i|+{1\over\varepsilon'}\sum_{i>J_{l'-1}}|u_i|b_i)\right]^2d\gamma(u)
\]
is finite by using inequality \eqref{eq:exp1} (see also \cite{hoang2020analysis} Proposition B4). As $h_0(u)<2^{-\min\{l,l'\}}$ for $u\in U_2^{ll'}$ and $\int_U(h_0(u))^{-2})d\gamma(u)$ is finite, $\gamma(U_2 ^{ll'})\le c2^{-2\min\{l,l'\}}$. We then have
\begin{align*}
&\int_U\exp(2c\sum_{i=1}^\infty|u_i|(b_i+\bar b_i)+{2\over \min\{\varepsilon,\varepsilon'\}}\sum_{i>\max\{J_{l-1},J_{l'-1}\}}|u_i|b_i){\cal I}_{U_2^{ll'}}(u)d\gamma(u)\\
&\le \left(\int_U\exp(4c\sum_{i=1}^\infty|u_i|(b_i+\bar b_i)+{4\over \min\{\varepsilon,\varepsilon'\}}\sum_{i>\max\{J_{l-1},J_{l'-1}\}}|u_i|b_i)d\gamma(u)\right)^{1/2}\gamma(U_2^{ll'})^{1/2}\le c2^{-\min\{l,l'\}}.
\end{align*}
Similarly, as $\int_U h_0(u)^{-4}d\gamma(u)$ is finite, $\gamma(U_3^{ll'})\le c2^{-4\max\{l,l'\}}$. Thus
\[
\int_U\exp\left(c\sum_{i=1}^\infty |u_i|b_i\right){\cal I}_{U_3^{ll'}}d\gamma(u)\le c2^{-2\max\{l,l'\}}\le c2^{-(l+l')}.
\]
Thus $\mathbb{E}^{\gamma}[|A_1^{ll'}|^2]d\gamma(u)\le c2^{-(l+l')}$. 
The proof of the convergence of the MLMCMC sampling procedure follows exactly from that in Appendix A of \cite{hoang2020analysis}. For the pCN sampler, if we assume the spectral gap result of Hairer et al. \cite{hairer2014spectral}, then the convergence rate of the MLMCMC sampling procedure holds.

\end{document}